\documentclass[envcountsame]{llncs}

\usepackage[utf8]{inputenc}

\usepackage{hyperref}
\usepackage{amsmath,amssymb,mathtools}
\usepackage[capitalise]{cleveref}

\usepackage{todonotes}
\usepackage{tikz}
\usepackage{todonotes}

\usepackage[linesnumbered, vlined, ruled]{algorithm2e}
\crefname{algocf}{alg.}{algs.}
\Crefname{algocf}{Algorithm}{Algorithms}

\usepackage{multirow}
\usepackage{booktabs}
\usepackage{tikz}
\usetikzlibrary{calc, positioning, arrows, arrows.meta, decorations.pathmorphing, decorations.pathreplacing, decorations.text, words, backgrounds, babel}

\newcommand{\rev}[1]{\mathpalette\harpoonvec{#1}}
\newcommand{\harpvecsign}{\scriptscriptstyle\leftharpoonup}
\newcommand{\harpoonvec}[2]{%
  \ifx\displaystyle#1\doalign{$\harpvecsign$}{#1#2}\fi
  \ifx\textstyle#1\doalign{$\harpvecsign$}{#1#2}\fi
  \ifx\scriptstyle#1\doalign{\scalebox{.6}[.9]{$\harpvecsign$}}{#1#2}\fi
  \ifx\scriptscriptstyle#1\doalign{\scalebox{.5}[.8]{$\harpvecsign$}}{#1#2}\fi
}
\newcommand{\doalign}[2]{%
 {\vbox{\offinterlineskip\ialign{\hfil##\hfil\cr#1\cr$#2$\cr}}}%
}

\def\N{\mathbb{N}}

\renewcommand\th[1]{{#1}\textsuperscript{th}}
\def\nth#1{\th{#1}}

\DeclareMathOperator\ar{ar}
\DeclareMathOperator\re{re}
\DeclareMathOperator\m{m}
\DeclareMathOperator\ra{\rev\ar}
\DeclareMathOperator\er{\rev\re}
\DeclareMathOperator\w{\rev\m}
\DeclareMathOperator\W{\rev{\mathrm M}}
\DeclareMathOperator\alphabet{alph}
\DeclareMathOperator\letters{alph}

\DeclareMathOperator\ScatFact{ScatFact}
\DeclareMathOperator\core{core}
\DeclarePairedDelimiter{\abs}{\lvert}{\rvert}
\DeclarePairedDelimiter{\norm}{\lVert}{\rVert}
\DeclarePairedDelimiter{\paren}{(}{)}

\def\ta{{\normalfont\texttt{a}}}
\def\tb{{\normalfont\texttt{b}}}
\def\tc{{\normalfont\texttt{c}}}
\def\td{{\normalfont\texttt{d}}}
\def\te{{\normalfont\texttt{e}}}

\def\tx{{\normalfont\texttt{x}}}
\def\ty{{\normalfont\texttt{y}}}

\definecolor{cau}{RGB}{156,10,125} 

\newif\ifpaper
\paperfalse 

\ifpaper
\else
\usepackage{bookmark}
\fi

\spnewtheorem{numclaim}{Claim}{\bfseries}{\rmfamily}
\crefname{numclaim}{claim}{claims}
\Crefname{numclaim}{Claim}{Claims}

\title{\texorpdfstring{\(\alpha\)-\(\beta\)}{αβ}-Factorization and the Binary Case of Simon's Congruence}
\author{Pamela Fleischmann\inst{1} \and Jonas Höfer\inst{2} \and Annika Huch\inst{1} \and Dirk Nowotka\inst{1}}

\authorrunning{P. Fleischmann, J. Höfer, A. Huch and D. Nowotka}

\institute{Kiel University, Germany, \email{$\{$fpa,dn$\}$@informatik.uni-kiel.de, stu216885@mail.uni-kiel.de}
\and University of Gothenburg, Sweden, \email{jonas.hofer@gu.se}}

\begin{document}
\maketitle
\begin{abstract}
  In 1991 Hébrard introduced a factorization of words that turned out to be a powerful tool for the investigation of a word's scattered factors (also known as (scattered) subwords or subsequences).
  Based on this, first Karandikar and Schnoebelen introduced the notion of $k$-richness and later on Barker et al. the notion of $k$-universality.
  In 2022 Fleischmann et al. presented at DCFS a generalization of the arch factorization by intersecting the arch factorization of a word and its reverse.
  While the authors merely used this factorization for the investigation of shortest absent scattered factors, in this work we investigate this new $\alpha$-$\beta$-factorization as such.
  We characterize the famous Simon congruence of $k$-universal words in terms of $1$-universal words.
  Moreover, we apply these results to binary words.
  In this special case, we obtain a full characterization of the classes and calculate the index of the congruence.
  Lastly, we start investigating the ternary case, present a full list of possibilities for $\alpha\beta\alpha$-factors, and characterize their congruence.
\end{abstract}

\section{Introduction}

A \emph{scattered factor}, \emph{subsequence}, \emph{subword} or \emph{scattered subword} of a word \(w\) is a word that is obtained by deleting any number of letters from \(w\) while preserving the order of the remaining letters.
For example, \(\texttt{oiaoi}\) and \(\texttt{cmbntrcs}\) are both scattered factors of \(\texttt{combinatorics}\).
In contrast to a factor, like \(\texttt{combinat}\), a scattered factor is not necessarily contiguous.
Note that a scattered factor \(v\) can occur in different ways inside a word \(w\), for example, \(\ta\tb\) occurs in \(\ta\ta\tb\) as \(\underline\ta\ta\underline\tb\) and \(\ta\underline\ta\underline{\tb}\) as marked by the lines below the letters.
The relation of \(u\) being a scattered factor of \(v\) is a partial order on words.

In this paper, we focus on the congruence relation \(\sim_k\) for \(k\in\mathbb N_0\) which is known as Simon's congruence \cite{simon1972thesis}.
For two words, we have \(u \sim_k v\) iff they share all scattered factors up to length \(k\).
Unions of the congruence classes of this relation are used to form the \emph{piecewise testable languages} (first studied by Simon~\cite{simon1972thesis}), which are a subclass of the regular languages (they are even subregular).

A long-standing open question, posed by Sakarovitch and Simon~\cite{lothaire1997:subwords}, is the exact structure of the congruence classes of \(\sim_k\) and the index of the congruence relation itself.
Two existing results include a characterization of the congruence in terms of a special common upper bound of two words~\cite[Lemma 6]{DBLP:conf/automata/Simon75}, as well as a characterization of the (not unique) shortest elements of the congruence 
classes~\cite[Theorem 6.2.9]{lothaire1997:subwords} and \cite{simon1972thesis,DBLP:conf/mfcs/FleischerK18,DBLP:conf/dlt/BarkerFHMN20}.
The index of the relation is described asymptotically by Karandikar et al.~\cite{DBLP:journals/ipl/KarandikarKS15}.
Currently, no exact formula is known.
One approach for studying scattered factors in words is based on the notion of \emph{scattered factor universality}.
A word \(w\) is called \emph{\(\ell\)-universal} if it contains all words of length \(\ell\) over a given alphabet as scattered factors. For instance, the word
$\mathtt{alfalfa}$\footnote{\href{https://en.wikipedia.org/wiki/Alfalfa}{Alfalfa (Medicago sativa)} is a plant whose name means \emph{horse food} in Old Persian.} is $2$-universal since it contains all words of length two over the alphabet $\{\ta,\mathtt{l},\mathtt{f}\}$ as scattered factors.
Barker et al. and Fleischmann et al.~\cite{DBLP:conf/dlt/BarkerFHMN20,fleischmann2021scattered} study the universality of words, as well as how the universality of a word changes when considering repetitions of a word.
Fleischmann et al.~\cite{DBLP:conf/dcfs/FleischmannHHMN22} investigate the classes of Simon's congruence separated by the number of shortest absent scattered factors, characterize the classes for arbitrary alphabets for some fixed numbers of shortest absent scattered factors and give explicit formulas for these subsets. The shortest absent scattered factors of $\mathtt{alfalfa}$ are $\mathtt{fff}$, $\mathtt{ffl}$ $\mathtt{lll}$, and $\mathtt{fll}$.
A main tool in this line of research is  a newly introduced factorization, known as the \emph{\(\alpha\)-\(\beta\)-factorization} \cite{DBLP:conf/dcfs/FleischmannHHMN22} which is based on the arch factorization by Hébrard~\cite{DBLP:journals/tcs/Hebrard91}. The arch factorization factorizes a word into factors of minimal length containing the complete alphabet. The $\alpha$-$\beta$-factorization takes also the arch factorization of the reversed word into account.
Kosche et al.~\cite{kosche2021absent} implicitly used this factorization to determine shortest absent scattered factors in words.
In this paper, we study this factorization from a purely combinatorial point of view.
The most common algorithmic problems regarding Simon's congruence are \textsc{SimK} (testing whether two words \(u,v\) are congruent for a fixed \(k\)) and \textsc{MaxSimK} (the optimization problem of finding the largest \(k\) such that they are congruent).
The former was approached by finding the (lexicographical least element of the) minimal elements of the congruence classes of \(u\) and \(v\).
Results regarding normal forms and the equation $pwq\sim_kr$ for given words $p,q,r$ can be found in \cite{pach,pach2}.
The computation of the normal form was improved first by Fleischer et al.~\cite{DBLP:conf/mfcs/FleischerK18} and later by Barker et al.~\cite{DBLP:conf/dlt/BarkerFHMN20}.
The latter was approached in the binary case by Hébrard~\cite{DBLP:journals/tcs/Hebrard91}, and was solved in linear time using a new approach by Gawrychowski et al.~\cite{DBLP:conf/stacs/GawrychowskiKKM21}. 
A new perspective on $\sim_k$ was recently given by Sungmin Kim et al. \cite{DBLP:conf/dcfs/KimHKS22,DBLP:conf/isaac/KimKH22} when investigating the congruence's closure and pattern matching w.r.t. $\sim_k$.

{\em Our Contribution.}
We investigate the \(\alpha\)-\(\beta\)-factorization as an object of independent interest and give necessary and sufficient conditions for the congruence of words in terms of their factors.
We characterize $\sim_k$ in terms of  \(1\)-universal words through their \(\alpha\beta\alpha\)-factors.
We use these results to characterize the congruence classes of binary words and their cardinality, as well as to calculate the index in this special case.
Moreover, we give a short and conceptually straightforward algorithm for \textsc{MaxSimK} for binary words.
Lastly, we start to transfer the previous results to the ternary alphabet.

{\em Structure of the Work.}
First, in \Cref{sec:preliminaries} we establish basic definitions and notation.
In \Cref{sec:alpha-beta}, we give our results regarding the \(\alpha\)-\(\beta\)-factorization for arbitrary alphabets, including the characterization of the congruence of words w.r.t. \(\sim_k\) in terms of their \(\alpha\beta\alpha\)-factors.
Second, in \Cref{sec:binary}, we present our results regarding binary words.
We characterize the congruence classes of binary words in terms of their \(\alpha\)- and \(\beta\)-factors, and apply them to calculate the index of \(\sim_k\) in this special case.
Third, in \Cref{sec:ternary}, we consider a ternary alphabet and investigate the cases for the \(\beta\)-factors.
Last, in \Cref{sec:conclusion}, we conclude and give ideas for further research.


\section{Preliminaries}\label{sec:preliminaries}

We set \(\mathbb N \coloneqq \{1,2,3,\dots\}\) and \(\mathbb N_0 \coloneqq \{0\} \cup \mathbb N\) as well as \([m]\coloneqq\{1,\ldots,m\}\) and \([m]_0\coloneqq\{0\}\cup [m]\).
We denote disjoint unions by \(\sqcup\).
If there exists a bijection between two sets \(A, B\), then we write \(A \cong B\).
An {\em alphabet} is a finite set $\Sigma$ whose elements are called {\em letters}.
An alphabet of cardinality $i\in\N$ is abbreviated by $\Sigma_i$.
A {\em word} $w$ is a finite sequence of letters from $\Sigma$ where $w[i]$ denotes the \nth{$i$} letter of $w$.
The set of all words over the alphabet $\Sigma$ is denoted by $\Sigma^{\ast}$ and the {\em empty word} by $\varepsilon$.
Set $\Sigma^+\coloneqq\Sigma^{\ast}\setminus\{\varepsilon\}$.
The {\em length} $|w|$ of $w$ is the number of letters in $w$, i.e., $|\varepsilon|=0$.
We denote the set of all words of length $k\in\N_0$ by $\Sigma^k$ and set $\Sigma^{\leq k}\coloneqq\{w\in\Sigma^{\ast}|\,|w|\leq k\}$.
Set $\letters(w)\coloneqq\{w[i] \in\Sigma \mid i\in[|w|]\}$.
Set $|w|_{\ta}\coloneqq|\{i\in[|w|]\mid w[i]=\ta\}|$ for all $\ta\in\Sigma$.
A word $u\in\Sigma^{\ast}$ is called {\em factor} of $w\in\Sigma^{\ast}$ if there exist $x,y\in\Sigma^{\ast}$ with $w=xuy$.
In the case that $x=\varepsilon$, $u$ is called {\em prefix} of $w$ and {\em suffix} if $y=\varepsilon$.
The factor of $w$ from its \nth{$i$} letter to its \nth{$j$} letter is denoted by $w[i..j]$ for $1\leq i\leq j\leq|w|$. For $j < i$ we define $w[i..j] \coloneqq \varepsilon$.
If $w=xy$ we write $x^{-1}w$ for $y$ and $wy^{-1}$ for $x$.
For $u\in\Sigma^{\ast}$ we set $u^0\coloneqq\varepsilon$ and inductively $u^{\ell}\coloneqq uu^{\ell-1}$ for all $\ell\in\N$.
For $w\in\Sigma^{\ast}$ define $w^R$ as $w[|w|]\cdots w[1]$.
For more background information on \emph{combinatorics on words} see \cite{lothaire1997}.

Now, we introduce the main notion of our work, the scattered factors also known as (scattered) subwords or subsequence (also cf. \cite{lothaire1997:subwords}).

\begin{definition}
  A word $u\in\Sigma^{\ast}$ of length $n\in\N_0$ is called a \emph{scattered factor} of $w\in\Sigma^{\ast}$ if there exist $v_0,\ldots, v_{n}\in\Sigma^{\ast}$ with $w=v_0u[1]v_1\cdots v_{n-1}u[n]v_{n}$, denoted by $u\preceq w$.
  Let $\ScatFact(w)\coloneqq \{v\in\Sigma^{\ast}\mid v\preceq w\}$ as well as $\ScatFact_k(w)\coloneqq \ScatFact(w)\cap\Sigma^k$ and $\ScatFact_{\leq k}(w)\coloneqq \ScatFact(w)\cap \Sigma^{\leq k}$.
\end{definition}

For instance, we have $\mathtt{and}\preceq\mathtt{agenda}$ but $\mathtt{nada}\not\preceq\mathtt{agenda}$. For comparing words w.r.t. their scattered factors, Simon introduced a congruence relation nowadays known as {\em Simon's congruence} \cite{simon1972thesis}. Two words are called Simon $k$-congruent, if they have the same set of scattered factors up to length $k$. We refer to this $k$ as the \emph{level} of the congruence. This set is the {\em full $k$-spectrum} of a word, whereas the {\em $k$-spectrum} only contains all scattered factors of exactly length $k$.

\begin{definition}
  Let $k\in\N$. Two words \(u,v\in\Sigma^{\ast}\) are called \emph{Simon \(k\)-congruent (\(u\sim_k v\))} iff \(\ScatFact_{\leq k}(u)=\ScatFact_{\leq k}(v)\). Let $[u]_{\sim_k}$ denote the congruence class of \(u\) w.r.t. \(\sim_k\).
\end{definition}

For instance, over $\Sigma=\{\ta,\tb\}$, the words $\ta\tb\ta\ta\tb\ta$ and $\tb\ta\ta\tb$ are Simon $2$-congruent since both contain each all words up to length $2$ as scattered factors. On the other hand, they are not Simon $3$-congruent since we have $\ta\ta\ta\preceq\ta\tb\ta\ta\tb\ta$ but $\ta\ta\ta\not\preceq\tb\ta\ta\tb$.

Starting in \cite{DBLP:journals/ipl/KarandikarKS15,journals/lmcs/KarandikarS19} and \cite{karandikar2016height} special $k$-spectra were investigated in the context of piecewise testable languages: the {\em rich} resp. {\em $k$-rich} words. This work was pursued from the perspective of the universality problem for languages in \cite{DBLP:conf/dlt/BarkerFHMN20,day2021edit,DBLP:conf/stacs/GawrychowskiKKM21,fleischmann2021scattered} with the new notion of {\em $k$-universal} words.

\begin{definition}
  A word $w\in\Sigma^{\ast}$ is called {\em $k$-universal} w.r.t. $\Sigma$ if $\ScatFact_k(w)=\Sigma^k$.
  The maximal $k$ such that $w$ is $k$-universal is denoted by $\iota_{\Sigma}(w)$ and called $w$'s \emph{universality index}.
\end{definition}

\begin{remark}
  If we investigate a single word $w\in\Sigma^{\ast}$, we assume $\Sigma=\letters(w)$ implicitly and omit the $\Sigma$ as index of $\iota$. 
\end{remark}

In \cite{DBLP:conf/dcfs/FleischmannHHMN22} the notion of universality was extended to $m$-nearly $k$-universal words, which are words where exactly $m$ scattered factors of length $k$ are absent, i.e., $|\ScatFact_k(w)|=|\Sigma|^k-m$. In the last section of their paper the authors introduce a factorization of words based on the arch factorization (cf. \cite{DBLP:journals/tcs/Hebrard91}) in order to characterize the $1$-nearly $k$-universal words with $\iota(w)=k-1$. This work is closely related to the algorithmic investigation of shortest absent scattered factors \cite{DBLP:conf/stacs/GawrychowskiKKM21}.
Therefore, we introduce first the arch factorization and based on this the $\alpha$-$\beta$-factorization from \cite{DBLP:conf/dcfs/FleischmannHHMN22}. An arch is a factor of minimal length (when read from left to right) containing the whole alphabet. Consider the word $w=\mathtt{abaccaabca}$. This leads to the arch factorization $(\ta\tb\ta\tc)\cdot (\tc\ta\ta\tb)\cdot\tc\ta$ where the arches are visualized by the brackets. 

\begin{definition}
  For a word \(w \in \Sigma^*\) the \emph{arch factorization} is given by \(w \eqqcolon \ar_1(w) \cdots \ar_k(w) \re(w)\) for \(k \in \mathbb{N}_0\) with
  \(\letters(\ar_i(w))=\Sigma\) for all \(i \in [k]\), the last letter of \(\ar_i(w)\) occurs exactly once in \(\ar_i(w)\) for all $i \in [k]$, and
\(\letters(\re(w)) \subset \Sigma\).
  The words \(\ar_i(w)\) are called \emph{arches} of \(w\) and \(\re(w)\) is the \emph{rest} of \(w\).
  Define the \emph{modus of \(w\)} as \(\m(w) \coloneqq \ar_1(w)[\vert\ar_1(w)\vert] \cdots\ar_k(w)[\vert\ar_k(w)\vert] \in \Sigma^k\).
  For abbreviation let $\ar_{i..j}(w)$ denote the concatenation
  from the \nth{$i$} arch to the \nth{$j$} arch.
\end{definition}

The following remark is a direct consequence of the combination of the \(k\)-universality and the arch factorization.

\begin{remark}\label{rem:simon-refines-arch-num-equiv}
  Let \(w,w'\in\Sigma^*\) such that \(w \sim_k w'\) for some \(k\in\mathbb N_0\).
  Then either both \(w,w'\) have \(k\) or more arches or they both have less than \(k\) and the same number of arches. Moreover, we have $\iota(w)=k$ iff $w$ has exactly $k$ arches.
\end{remark}

A generalization of the arch factorization was introduced in \cite{DBLP:conf/dcfs/FleischmannHHMN22} inspired by \cite{kosche2021absent}.
In this factorization not only the arch factorization of a word $w$ but also the one of $w^R$ is taken into consideration.
If both arch factorisations, i.e., the one of $w$ and the one of $w^R$ are considered simultaneously, we get overlaps of the arches and special parts which start at a modus letter of a {\em reverse} arch and end in a modus letter of an arch. 
For better readability, we use a specific notation for the arch factorisation of $w^R$ where we read the parts from left to right: let 
\(\ra_i(w)\coloneqq(\ar_{\iota(w)-i+1}(w^R))^R\) the \nth{$i$} {\em reverse arch}, let \(\er(w) \coloneqq (\re(w^R))^R\) the {\em reverse rest}, and define the {\em reverse modus} $\w(w)$ as $\m(w^R)^R$.

\begin{figure}[tbh]
	\centering
	\begin{tikzpicture}[scale=0.8]
    \draw[point] (11, .5) -- ++(0, -.4);
    \draw[uword] (11, .5) -- ++(2, 0) node {\(\re(w)\)};
    \draw[point] (13, .5) -- ++(0, -.4);
    \foreach \x/\w/\i in {0/3/1,3/2/2,5/4/3,9/2/4}{
    	\draw[uword] (\x,0.5) -- ++(\w, 0) node {\(\ar_{\i}(w)\)};
    	\draw[point] (\x, .5) -- ++(0, -.4);
    }
    \draw[point] (2, -.5) -- ++(0, .4);
    \draw[point] (0, -.5) -- ++(0, .4);
    \draw[dword] (0, -.5) -- ++(2, 0) node {\(\er(w)\)};
    \foreach \x/\w/\i in {2/2/1,4/3/2,7/3/3,10/3/4}{
    	\draw[dword] (\x,-0.5) -- ++(\w, 0) node {\(\ra_{\i}(w)\)};
    	\draw[point] (\x+\w, -.5) -- ++(0, .4);
    }
    \draw[word] (0,0) -- ++(2, 0) node [innerLabel] {\(\alpha_{0}\)} node [pos=0, left] {\(w\)};
    \foreach \x/\w/\i in {3/1/1,5/2/2,9/1/3,11/2/4}{
    	\draw[word] (\x,0) -- ++(\w, 0) node [innerLabel] {\(\alpha_{\i}\)};
    }
    \foreach \x/\w/\i in {2/1/1,4/1/2,7/2/3,10/1/4}{
    	\draw[word] (\x,0) -- ++(\w, 0) node [innerLabel] {\(\beta_{\i}\)};
    }
	\end{tikzpicture}
	\caption[\(\alpha\)-\(\beta\)-Factorization]{\(\alpha\)-\(\beta\)-Factorization of a word \(w\) with \(4\) arches.}
	\label{fig:alpha-beta-factorization}
\end{figure}
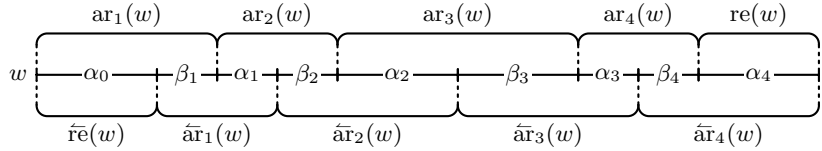

\begin{definition}

  For \(w\in\Sigma^{\ast}\) define \(w\)'s \emph{\(\alpha\)-\(\beta\)-factorization} (cf. \Cref{fig:alpha-beta-factorization}) by \(w\eqqcolon\alpha_0\beta_1\alpha_1\cdots \alpha_{\iota(w)-1}\beta_{\iota(w)}\alpha_{\iota(w)}\) with 
  \(\ar_i(w)=\alpha_{i-1}\beta_i\) and \(\ra_i(w)=\beta_i\alpha_i\) for all $i\in[\iota(w)]$, \(\er(w)=\alpha_0\), as well as $\re(w)=\alpha_{\iota(w)}$.
    Define $\core_i\coloneqq \varepsilon$ if $|\beta_i|\in\{1,2\}$ and $\core_i=\beta_i[2..|\beta_i|-1]$ otherwise, i.e., as the $\beta_i$ without the associated letters of the modus and reverse modus.
\end{definition}
For example, consider $w= \mathtt{bakebananacake} \in \{\mathtt{a,b,c,k,e}\}^*$. We get $\ar_1(w) = \mathtt{bakebananac}$, $\re(w)= \mathtt{ake}$ and $\ra_1(w) = \mathtt{bananacake}$, $\er(w) = \mathtt{bake}$.
Thus, we have $\alpha_{0} = \mathtt{bake}$, $\beta_1 = \mathtt{bananac}$ and $\alpha_1 = \mathtt{ake}$. Moreover, we have $\m(w)=\tc$ and $\w(w)=\tb$. This leads to $\core_1 = \mathtt{anana}$.

\begin{remark}
In contrast to the arch factorization, the \(\alpha\)-\(\beta\)-factorization is left-right-symmetric. 
  Note that the \nth{\(i\)} reverse arch always starts inside the \nth{\(i\)}  arch since otherwise an arch or the rest would contain at least two reverse arches or a complete arch and thus the arch would contain the complete alphabet more than once or once.
\end{remark}  
  
 For better readability, we do not parametrize the $\alpha_i$ and $\beta_i$ by $w$. Instead, we denote the factors according to the word's name, i.e. $\tilde{\alpha}_i\tilde{\beta}_{i+1}$ is an arch of $\tilde{w}$.

 \begin{remark}  
 Since $|\letters(\alpha_i)|\leq|\Sigma|-1$ we can build the arch factorization of $\alpha_i$ w.r.t.~some \(\Omega\) with \(\letters(\alpha_i) \subseteq \Omega \in \binom{\Sigma}{\vert \Sigma \vert -1}\).
  This yields the same factorization for all $\Omega$ because either \(\letters(\alpha_i) = \Omega\) or \(\letters(\alpha_i) \subset \Omega\) and thus \(\re(\alpha_i) = \alpha_i\).
\end{remark}

Last, we recall three lemmata regarding Simon's congruence which we need for our results.
The first lemma shows that if we prepend or append a sufficiently universal word to two congruent words each, we obtain congruent words with an increased level of congruence.

\begin{lemma}[{\cite[Lemma 4.1]{DBLP:journals/ipl/KarandikarKS15}\cite[Lemma 3.5]{journals/lmcs/KarandikarS19}}]\label{prop:extend-equiv-words-by-univ-word-increases-equiv}
  Let \(w, \tilde w \in\Sigma^*\) such that \(w \sim_k \tilde w\), then for all \(u, v\in \Sigma^*\) we have \(uwv \sim_{\iota(u) + k + \iota(v)} u\tilde wv\).
\end{lemma}

The next lemma characterizes the omittance of suffixes when considering words up to \(\sim_k\).

\begin{lemma}[{\cite[Lemma 3]{DBLP:conf/automata/Simon75}}]\label{prop:simon:suffix-elim}
  Let \(u,v\in\Sigma^+\) and \(\tx\in\Sigma\).
  Then, \(uv \sim_k u\) iff there exists a factorization \(u = u_1u_2\cdots u_k\) such that \(\alphabet(u_1) \supseteq \alphabet(u_2) \supseteq \ldots \supseteq \alphabet(u_k) \supseteq \alphabet(v)\).
\end{lemma}

The last lemma characterizes letters which can be omitted when we consider words up to \(\sim_k\).
The last two of its conditions follow from the previous lemma.

\begin{lemma}[{\cite[Lemma 4]{DBLP:conf/automata/Simon75}}]\label{prop:simon:letter-elim}
  Let \(u,v\in\Sigma^*\) and \(\tx\in\Sigma\).
  Then, \(uv\sim_k u\tx v\) iff there exist \(p,p'\in\mathbb N_0\) with \(p + p' \geq k\) and \(u\tx \sim_p u\) and \(\tx v\sim_{p'} v\).
\end{lemma}

\section{\texorpdfstring{\(\alpha\)-\(\beta\)}{αβ}-Factorization}\label{sec:alpha-beta}
In this section, we investigate the $\alpha$-$\beta$-factorization based on results of \cite{DBLP:journals/ipl/KarandikarKS15} in the relatively new light of factorizing an arch into an $\alpha$ and a $\beta$ part. The main result states that it suffices to look at $1$-universal words in order to gain the information about the congruence classes of $\sim_k$.

\begin{remark}
By the left-right symmetry of the $\alpha$-$\beta$-factorisation, it suffices to prove most of the claims only for one direction (reading the word from left to right) and the other direction (reading the word from right to left) follows immediately. Thus, these claims are only given for one direction and it is not always mentioned explicitly that the analogous claim holds for the other direction.
\end{remark}

Our first lemma shows that {\em cutting of} $\ell$ arches from two $k$-congruent words each, leads to $(k-\ell)$-congruence. Here, we use the $\alpha$-$\beta$-factorization's symmetry.

\begin{lemma}\label{prop:simon-equiv-implies-equiv-aba-factors}
    Let \(w,\tilde w \in \Sigma^*\) with \(w \sim_k \tilde w\) and \(\iota(w) = \iota(\tilde w)<k\). Then we have  \(\alpha_i\beta_{i+1}\alpha_{i+1}\cdots \alpha_j \sim_{k - \iota(w) + j - i} \tilde\alpha_i\tilde\beta_{i+1}\tilde\alpha_{i+1}\cdots \tilde\alpha_j\) for all \(0 \leq i \leq j \leq \iota(w)\).
\end{lemma}

\ifpaper
\else 
\begin{proof}
  By symmetry, it suffices to show one inclusion of the \((k-1)\)-spectra.
  Let \(v \preceq \ar_1^{-1}(w)\cdot w = \ar_{2..\iota(w)}(w)\re(w)\) with \(|v| \leq k - 1\) and \(\texttt{m}_{\texttt{1}} \coloneqq \m(\tilde w)[1]\).
  Then \(\texttt{m}_{\texttt{1}} \cdot v \preceq w\) because \(\texttt{m}_{\texttt{1}} \preceq \ar_1(w)\).
  Since \(w\sim_k \tilde w\), we have \(\texttt{m}_{\texttt{1}} \cdot v \preceq \tilde w\) and therefore \(v \preceq \ar_1^{-1}(\tilde w) \cdot \tilde w\). This proves the first claim.
  
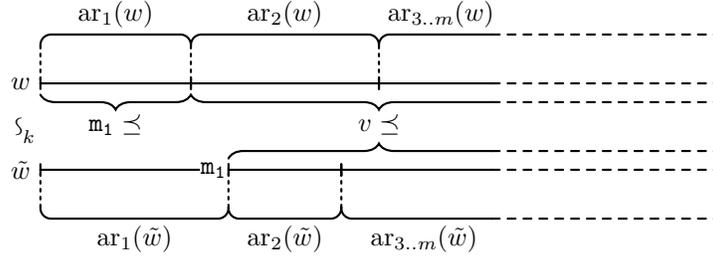
\begin{center}
    \begin{tikzpicture}
      \draw[word] (0, 0) -- ++(2, 0) node[outerLabelStart] (l1) {\(w\)};
      \draw[word] (0, -1.15) -- ++(2.5, 0) node[outerLabelStart] (l2) {\(\tilde w\)} node [innerLabelEnd] {\texttt{m\textsubscript{1}}};
  
      \node[rotate=90,xscale=-1, inner sep=1pt] (t) at ($(l1)!.5!(l2)$) {\(\sim\)};
      \node[font=\scriptsize,yshift=.5mm] at (t.south east) {\(k\)};
  
      \draw[uword] (0, .5) -- ++(2, 0) node {\(\ar_1(w)\)};
      \draw[dword] (0, -1.65) -- ++(2.5, 0) node {\(\ar_1(\tilde w)\)};
  
      \draw[word] (2, 0) -- ++(2.5, 0) node {};
      \draw[word] (2.5, -1.15) -- ++(1.5, 0) node {};
  
      \draw[uword] (2, .5) -- ++(2.5, 0) node {\(\ar_2(w)\)};
      \draw[dword] (2.5, -1.65) -- ++(1.5, 0) node {\(\ar_2(\tilde w)\)};
  
      \draw[brace] (2, -.15) -- ++(-2, 0) node [below] {\(\texttt{m\textsubscript{1}} \preceq \)};
      \begin{scope}
        \clip[overlay] (0, -2.5) rectangle (6, 1.5);
        \draw[brace] (7, -.15) -- ++(-5, 0) node [below] {\(v\preceq\)};
        \draw[brace] (2.5, -1) -- ++(4, 0) node {};
  
        \draw[word] (4.5, 0) -- ++(2, 0);
        \draw[word] (4, -1.15) -- ++(2.2, 0);
  
        \draw[uword] (4.5, .5) -- ++(1.65, 0) node {\(\ar_{3..m}(w)\)};
        \draw[dword] (4, -1.65) -- ++(2.2, 0) node {\(\ar_{3..m}(\tilde w)\)};
      \end{scope}
  
      \foreach \x in {0,2.5,4}{ \draw[point] (\x, -1.65) -- ++(0, 0.475); }
      \foreach \x in {0,2,4.5}{ \draw[point] (\x, 0.5) -- ++(0, -0.475); }
  
      \foreach \y in {.65,0,-0.25,-0.9,-1.15,-1.8}{ \draw [dashed, line cap=round, thick] (6, \y) -- ++(3, 0); }
    \end{tikzpicture}
  \end{center}
  The second claim follows by left-right-symmetry and induction by repeatedly cutting \(\alpha\beta\)-pairs from the left and \(\beta\alpha\)-pairs form the right.\qed
\end{proof}

\fi

The following proposition shows that two words having exactly the same $\beta$-factors are $k$-congruent iff the corresponding \(\alpha\)-factors are congruent at a smaller level.
The proof uses a similar idea to the one presented by Karandikar et al.~\cite[Lemma 4.2]{DBLP:journals/ipl/KarandikarKS15}.

\begin{proposition}\label{prop:beta-equal-implies-equiv-iff-alpha-equiv}
  For all \(w,\tilde w\in\Sigma^*\) with \(m \coloneqq \iota(w) = \iota(\tilde w) < k\) such that \(\beta_i = \tilde \beta_i\) for all \(i\in[m]\), we have \(w \sim_k \tilde{w}\) iff \(\alpha_i \sim_{k - m} \tilde\alpha_i\) for all \(i\in[m]_0\).
\end{proposition}
\ifpaper
\else 
\begin{proof}
  By Lemma~\ref{prop:simon-equiv-implies-equiv-aba-factors}, we directly obtain one direction.
  Therefore, let \(w,\tilde w\in\Sigma^*\) such that \(\beta_i = \tilde\beta_i\) for all \(i\in[m]\) and \(\alpha_i \sim_{k - m} \tilde\alpha_i\) for all \(i\in[m]_0\).
  We obtain by Lemma~\ref{prop:extend-equiv-words-by-univ-word-increases-equiv} 
  \begin{align*}
    \tilde \alpha_0\beta_1\tilde \alpha_1\cdots \tilde \alpha_{i-1} \beta_{i}\alpha_{i}\beta_{i+1}\cdots \beta_{\iota(w)}\alpha_{\iota(w)}
    \sim_k
    \tilde \alpha_0\beta_1\tilde \alpha_1\cdots \tilde \alpha_{i-1} \beta_{i}\tilde \alpha_{i}\beta_{i+1}\cdots \beta_{\iota(w)}\alpha_{\iota(w)}
  \end{align*}
  for all $i\in[m]_0$.
  Thus, by transitivity of \(\sim_k\), we have \(w\sim_k \tilde{w}\).\qed
\end{proof}

\fi

As an immediate corollary, we obtain the following statement which allows us to normalize the \(\alpha\)-factors when proving congruence of words.
\begin{corollary}\label{prop:simon-equiv-iff-equiv-with-same-alpha}
  Let \(w,\tilde w\in\Sigma^*\) with \(m \coloneqq \iota(w) = \iota(\tilde w) < k\), then \(w \sim_k \tilde w\) iff \(\alpha_i \sim_{k-m} \tilde \alpha_{i}\) for all \(i\in[m]_0\) and for \(w' \coloneqq \alpha_0\tilde\beta_1\alpha_1\cdots \tilde\beta_m\alpha_m \) we have \(w\sim_k w'\).
\end{corollary}
\ifpaper
\else
\begin{proof}
  Note that, \(w'\) is in \(\alpha\)-\(\beta\)-factorization because the exchanged \(\alpha_i\) are equivalent, and thus they have the same alphabets.
  Assume \(w \sim_k \tilde w\).
  By \Cref{prop:simon-equiv-implies-equiv-aba-factors} we have \(\alpha_i \sim_{k-m} \tilde\alpha_{i}\).
  Define \(w'\) as above.
  We have \(w' \sim_k \tilde w\) by \Cref{prop:beta-equal-implies-equiv-iff-alpha-equiv} and therefore, by transitivity \(w \sim_k w'\).

  Now assume the converse.
  By \(\alpha_i' = \alpha_i \sim_{k-m} \tilde\alpha_{i}\) and \(\beta_j' = \tilde\beta_j\) for all \(i\in[m]_0, j\in[m]\), we obtain by \Cref{prop:beta-equal-implies-equiv-iff-alpha-equiv} that \(\tilde w \sim_k w'\).
  By the assumption and transitivity, we obtain \(w\sim_k \tilde w\). \qed
\end{proof}

\fi
Next, we show the central result for this section.
We can characterize the congruence of words by the congruence of their \(\alpha\beta\alpha\)-factors.
Therefore, it suffices to consider \(1\)-universal words in general.
Again, the proof uses \Cref{prop:extend-equiv-words-by-univ-word-increases-equiv} and is inspired by Karandikar et al.~\cite[Lemma 4.2]{DBLP:journals/ipl/KarandikarKS15} and repeatedly exchanges factors up to \(k\)-Simon congruence.

\begin{theorem}\label{prop:equiv-iff-all-aba-equiv}
  Let \(w, \tilde w\in\Sigma^*\) with \(m \coloneqq \iota(w) = \iota(\tilde w) < k\).
  Then, \(w \sim_k \tilde w\) iff \(\alpha_{i-1}\beta_i\alpha_i \sim_{k - m + 1} \tilde\alpha_{i-1}\tilde\beta_i\tilde\alpha_i\) for all \(i\in[m]\).
\end{theorem}

\ifpaper
\else
\begin{proof}
  Assume \(w \sim_k \tilde w\), then the congruences follow directly by \Cref{prop:simon-equiv-implies-equiv-aba-factors} for \(i,j\in\mathbb{N}_0\) with \(|j - i| = 1\).

  Assume \(\alpha_{i-1}\beta_i\alpha_i \sim_{k-m+1} \tilde\alpha_{i-1}\tilde\beta_{i}\tilde\alpha_{i}\) for all \(i\in[m]\).
  By \Cref{prop:simon-equiv-implies-equiv-aba-factors}, we obtain that \(\alpha_i \sim_{k - m} \tilde\alpha_i\) for all \(i \in [m]_0\).
  By \Cref{prop:simon-equiv-iff-equiv-with-same-alpha}, we have \(\alpha_{i-1}\beta_i\alpha_i \sim_{k-m+1} \alpha_{i-1}\tilde\beta_i\alpha_i\) for all \(i\in[m]\), and it suffices to show that \(w \sim_k \alpha_0\tilde\beta_1\alpha_1\cdots \beta_m\alpha_m\).
  Now, we have by repeated applications of \Cref{prop:extend-equiv-words-by-univ-word-increases-equiv} that
  \begin{equation}
    \alpha_0\beta_1\alpha_1\cdot\beta_2\alpha_2 \cdots \beta_m\alpha_m \sim_k \alpha_0\tilde\beta_1\alpha_1\cdot\beta_2\alpha_2 \cdots \beta_m\alpha_m \sim_k \ldots \sim_k \tilde w. \tag*{\qed}
  \end{equation}
\end{proof}

\fi
In the light of \Cref{prop:equiv-iff-all-aba-equiv}, in the following, we consider some special cases of these triples w.r.t. the alphabet of both involved $\alpha$. Hence, let \(w,\tilde w\in\Sigma^*\) with \(1 = \iota(w) = \iota(\tilde w)\).

\begin{proposition}\label{prop:both-alpha-empty}
  Let \(\alpha_0 = \alpha_1 = \tilde\alpha_0 = \tilde\alpha_1 = \varepsilon\).
  Then \(w \sim_k \tilde w\) iff \(k = 1\) or \(k \geq 2\), \(\m(w) = \m(\tilde w)\), \(\w(w) = \w(\tilde w)\), and \(\core_1 \sim_k \widetilde\core_1\).
\end{proposition}
\ifpaper
\else
\begin{proof}
  If \(k = 1\), the claim follows directly.
  If \(k \geq 2\), \(\m(w) = \m(\tilde w)\), \(\w(w) = \w(\tilde w)\) and \(w[2..|w|-1] = \core_1 \sim_k \widetilde\core_1 = \tilde w[2..|\tilde w|-1]\), then \(w \sim_k \tilde w\) follows directly from the fact that \(\sim_k\) is a congruence.

  Assume \(w \sim_k \tilde w\) and \(k \geq 2\).
  Suppose \(\m(w) \ne \m(\tilde w)\).
  Then \(\m(\tilde w)\m(w) \preceq w\) because \(\alphabet(w) = \Sigma\).
  By $w\sim_k\tilde{w}$, we get \(\m(\tilde w)\m(w) \preceq \tilde w = \ar_1(\tilde w)\), a contradiction.
  Therefore, \(\m(w) = \m(\tilde w)\) and by symmetry \(\w(w) = \w(\tilde w)\).
  Now the claim follows because \(\core_1 = \pi_\Omega(w) \sim_k \pi_\Omega(\tilde w) = \widetilde \core_1\) where \(\Omega = \Sigma\setminus\{\m(w),\w(w)\} = \Sigma\setminus\{\m(\tilde w),\w(\tilde w)\}\).\qed
\end{proof}

\fi

\begin{proposition}\label{prop:alpha-miss-same-single}
  Let  \(\letters(\alpha_0) = \letters(\alpha_1) = \letters(\tilde\alpha_0) = \letters(\tilde\alpha_1) \in \binom{\Sigma}{|\Sigma|-1}\).
  We have \(w \sim_k \tilde w\) iff \(\alpha_i \sim_{k-1} \tilde\alpha_i\) for all \(i\in[1]_0\).
\end{proposition}
\ifpaper
\else
\begin{proof}
  In both directions, we have \(\beta = \tilde\beta \in \Sigma\) by uniqueness of their first and last letter (which are identical).
  The claim follows from \Cref{prop:beta-equal-implies-equiv-iff-alpha-equiv}.\qed
\end{proof}

\fi

In the last two propositions, we considered special cases of congruence classes, where all words in such a congruence class have not only the same modus but also the same reverse modus. This is not necessarily always the case witnessed by 
$w=\ta \tb \ta \tb \te \ta \tb \ta \tb \cdot \ta\tb\te \tc\td \cdot \tc\td\tc\td\tc\td \sim_4 \ta\tb\ta\tb\te\ta\tb\ta\tb \cdot \tb\ta\te\td\tc \cdot \tc\td\tc\td\tc\td=\tilde{w}$ with $\m(w)=\td\neq \tc=\m(\tilde{w})$ and $\w(w)=\ta\neq \tb=\w(\tilde{w})$. This case occurs if one of the \(\alpha\)  satisfies $\alpha_0 \tx\sim_{k-1}\alpha_0$  and  the alphabet of \(\alpha_1\) factor is missing at least $\tx$ for all \(\tx \in \{ \w(\tilde w) \mid \tilde w \in [w]_{\sim_k} \}\).
The conditions for $\m(w)$ are analogous.
In the last proposition of this section, we show a necessary condition for the \(\alpha\)-factors of words which are congruent to words with a different modus.
The proof uses the same factorization as the proof of \Cref{prop:simon:suffix-elim} (cf. \cite[Lemma 3]{DBLP:conf/automata/Simon75}).
By identifying permutable factors, similar ideas also appear when characterizing the shortest elements in congruence classes (cf. \cite[Theorem 6.2.9]{lothaire1997:subwords}\cite[Proposition 6]{DBLP:conf/mfcs/FleischerK18}).

\begin{proposition}\label{prop:different-modi-force-alpha-univ-k-1}
  Let \(w\in \Sigma^*\) with \(\iota(w) = 1\), \(k \in\mathbb N\), and \(\W \coloneqq \{ \w(\tilde w)[1] \mid \tilde w \in [w]_{\sim_k} \}\), i.e., we capture all modus letters of words which are $k$-congruent to $w$.
  If \(|\W| \geq 2\) then there exists a factorization \(\alpha_0 \eqqcolon u_1\cdots u_{k-1}\) with \(\letters(u_1) \supseteq \ldots \supseteq \letters(u_{k-1}) \supseteq \W\).
\end{proposition}
\ifpaper
\else
\begin{proof}
  Let \(w, \W\) be as above and assume \(|\W| \geq 2\).
  We define the following factorization inductively.
  Let \(z \in \Sigma^*\).
  If \(\W\not\subseteq\letters(z)\), the factorization of \(z\) is \(v_1\coloneqq z\).
  Otherwise, set \(v_1\) as the shortest prefix of \(z\) containing \(\W\).
  Note that the last letter of \(v_1\) is unique.
  By induction on the length \(z\), there exists a factorization \(v_2\cdots v_n\) of \(v_1^{-1}z\).
  The factorisation of \(z\) is given by \(v_1v_2\cdots v_n\).
  Note that the last letters of \(v_1,\ldots,v_{n-1}\) are unique and \(\W \not\subseteq \letters(v_n)\).

  By factorizing \(\alpha_0\) this way, we obtain a factorization \(v_1v_2\cdots v_n\).
  Note that we have \(\letters(v_i) \supseteq \W\) exactly for \(i\in[n-1]\).
  If \(n \geq k\), define \(u_i \coloneqq v_i\) for \(i\in[k-2]\) and \(u_{k-1} \coloneqq v_{k-1}\cdots v_n\), and the claim follows.
  
  Therefore, suppose \(n < k\).
  Define \(m \coloneqq v_1[|v_1|]\cdots v_{n-1}[|v_{n-1}|] \in \Sigma^{\leq k - 2}\) as the unique last letters of \(v_1,\ldots,v_{n-1}\).
  Choose \(\tx,\ty\in \W\) distinct.
  Then there exist \(w_\tx, w_\ty \in [w]_{\sim_k}\) with \(\w(w_\tx) = \tx, \w(w_\ty) = \ty\) and assume by \Cref{prop:beta-equal-implies-equiv-iff-alpha-equiv} that they have the same \(\alpha\)-factors as \(w\).
  By definition, \(m \tx\ty \preceq w_\tx\).
  Since \(w_\tx \sim_k w_\ty\), we have \(m\tx\ty\preceq w_\ty\) and thus \(m\tx \preceq \alpha_0 = v_1\cdots v_n\) and thus \(\tx \preceq v_n\).
  Because \(\ty,\tx\) were chosen arbitrarily, we have \(\W \subseteq \letters(v_n)\), a contradiction against the construction of \(v_n\).
  This implies \(n \geq k\). \qed
\end{proof}

\fi

\section{The Binary Case of Simon's Congruence}\label{sec:binary}
In this section, we apply our previous results to the special case of the binary alphabet. Here, for $x\in\Sigma$, let $\overline{x}$ be the well defined other letter of $\Sigma$.
First, we characterize the congruence of binary words in terms of \(\alpha\)- and \(\beta\)-factors.
We show that in this scenario in each congruence class of a word $w$ with at most $k$ arches, we have $|\{ \w(\tilde w) \mid \tilde w \in [w]_{\sim_k} \}|=1$ (cf. Proposition~\ref{prop:different-modi-force-alpha-univ-k-1}).
We present results such that a full characterization of the structure of the classes in the binary case is given, implying as a byproduct a simple algorithm for \textsc{MaxSimK} in this special case (cf. \cite{DBLP:journals/tcs/Hebrard91}). Moreover, we can calculate \(|\Sigma_2^*\mathord{/}{\sim_k}|\).

\begin{proposition}\label{prop:binary:beta-alpha-cases}
  For all \(w\in\Sigma_2^*\), we have for all \(i\in[\iota(w)]\)\\
  1.  \(\beta_i \in \{\ta,\tb,\ta \tb,\tb \ta\}\),\\
  2.  if \(\beta_i = \tx\), then \(\alpha_{i-1}, \alpha_i \in \overline \tx^+\) with $\tx \in \Sigma_2$,\\
  3. if \(\beta_i = \tx\overline \tx\), then \(\alpha_{i-1} \in \tx^*\) and \(\alpha_i \in \overline \tx^*\) with $\tx \in \Sigma_2$.
  \end{proposition}

\ifpaper
\else
\begin{proof}
  By definition, the first and last letter of \(\beta_i\) are unique.
  Therefore, if \(\beta_i[1] = \beta_i[|\beta_i|]\) we have \(|\beta_i| = 1\).
  Furthermore, if \(\beta_i[1] \ne \beta_i[|\beta_i|]\) we have \(|\beta_i| = 2\)  because \(|\Sigma_2| = 2\).

  Because \(\alphabet(\alpha_i) \subset \Sigma_2\), the \(\alpha_i\) are unary words for all \(i\in[m]_0\).
  By symmetry, we only have to show the claim for \(\alpha_{i-1}\).
  The restrictions on the alphabet of \(\alpha_{i-1}\) follow directly from \(\ar_i(w) = \alpha_{i-1}\beta_i\) and the uniqueness of the last letter.
  Furthermore, if \(\beta_i = \tx \in \Sigma_2\) then \(\alpha_{i-1} \ne \varepsilon\) because \(\tx,\overline \tx\preceq \ar_i(w) = \alpha_{i-1}\beta_i\). \qed
\end{proof}

\fi

Thus, we get immediately that the $\alpha\beta\alpha$-factors are of the following forms: $\ta^{\ell_1+1}\tb\ta^{\ell_2+1}$, $\tb^{\ell_1+1}\ta\tb^{\ell_2+1}$, $\tb^{\ell_3}\tb\ta\ta^{\ell_4}$, or
$\ta^{\ell_3}\ta\tb\tb^{\ell_4}$ for some $\ell_1,\ell_2,\ell_3,\ell_4\in\N_0$. The following lemma shows that in the binary case the $k$-congruence of two words with identical universality less than $k$ leads to the same modi and same $\beta$.

\begin{lemma}\label{prop:binary:simon-equiv-implies-beta-equiv}
   Let \(w,w'\in\Sigma_2^*\) with \(w \sim_k w'\) and \(m \coloneqq \iota(w) = \iota(w') < k\), then \(\m(w) = \m(w')\) and thus, \(\beta_i = \beta_i'\) for all \(i\in[m]\).
\end{lemma}
\ifpaper
\else
\begin{proof}
  By \Cref{prop:equiv-iff-all-aba-equiv}, without loss of generality, we only consider \(m = 1 < k\).
  Furthermore, it suffices to show that \(\m(w)=\m(w')\) because then symmetry implies \(\w(w)=\w(w')\), and \(\m(w)\) and \(\w(w)\) fully determine the \(\beta\) factors by \Cref{prop:binary:beta-alpha-cases}.

  We know that the \(\alpha\)-factors are unary or empty by \Cref{prop:binary:beta-alpha-cases}.
  \Cref{prop:simon-equiv-implies-equiv-aba-factors} implies \(\alphabet(\alpha_i) = \alphabet(\alpha_i')\) for $i\in[1]_0$ since we have $w\sim_kw'$.
  If \(\alphabet(\alpha_0) \ne \emptyset\), then the arches \(\alpha_0\beta_1\) and \(\alpha_0'\beta_1'\) start with the same letter and thus \(\m(w)=\m(w')\).

  For the second case, assume \(\alpha_0 = \alpha_0' = \varepsilon\).
  If \(\beta_1 = \beta_1'\), then in particular \(\m(w) = \m(w')\), and the claim follows.
  Thus, suppose \(\beta_1 = \ty\overline \ty\) and \(\beta_1' = \overline \ty\ty\) for some \(\ty\in\Sigma_2\) by \Cref{prop:binary:beta-alpha-cases}.
  Thus, \(\alpha_1 \in \overline\ty^*\) and \(\alpha_1' \in \ty^*\).
  Hence, \(\alpha_1 = \alpha_1' = \varepsilon\) because \(\ty \notin \letters(\alpha_1) = \letters(\alpha_1') \not\ni \overline \ty\).
  This contradicts \(w\sim_k w'\) because \(w = \ty\overline\ty \not\sim_2 \overline\ty\ty = w'\). \qed
\end{proof}

\fi

Combining the Lemmata \ref{prop:binary:simon-equiv-implies-beta-equiv},~\ref{prop:simon-equiv-implies-equiv-aba-factors} and Proposition~\ref{prop:beta-equal-implies-equiv-iff-alpha-equiv} yields the following characterization of $\sim_k$ for binary words in terms of unary words and factors.

\begin{theorem}\label{prop:binary:simon-characterization}
  Let \(w,w'\in\Sigma_2^*\) such that \(m \coloneqq \iota(w) = \iota(w') < k\), then \(w\sim_k w'\) iff \(\beta_i = \beta_i'\) for all \(i\in[m]\) and \(\alpha_i \sim_{k - m} \alpha_i'\) for all \(i\in[m]_0\).
\end{theorem}

Using the characterization, we can also give an \(\mathcal O(|u|+|v|)\)-time algorithm for finding the largest \(k\) with \(u \sim_k v\) for \(u,v\in\Sigma_2^*\).
This special case was originally solved by Hébrard~\cite{DBLP:journals/tcs/Hebrard91} just considering arches.
Recently, a linear time algorithm for arbitrary alphabets was presented by Gawrychowski et al.~\cite{DBLP:conf/stacs/GawrychowskiKKM21}.
Nonetheless, we give \Cref{algo:bin-sim-k}, as it is a conceptually simple algorithm exploiting  that \(\alpha_i\) factors can be treated similar to \(\re(w)\) in the arch factorization.

\begin{algorithm}[hb]
  \caption{\textsc{MaxSimK} for binary words}
  \label{algo:bin-sim-k}
  \SetVlineSkip{0pt}
  \SetFuncSty{textsc}
  \SetKwFunction{ABFact}{\(\alpha\)-\(\beta\)-Fact}
  \KwIn{\(u, \tilde u\in\Sigma_2^*\)}
  \KwResult{if \(u = u'\) then \(\infty\) and otherwise the maximum \(k\) such that \(u \sim_k \tilde u\)}
  \((\alpha_0,\beta_1,\dots,\alpha_{\iota(u)}) \coloneqq \ABFact{\(u\)}\)\tcp*{w.r.t. \(\Sigma_2\)}
  \((\tilde\alpha_0,\tilde\beta_1,\dots,{\makebox[0pt][l]{\(\tilde\alpha_{\iota(\tilde u)}\)}\phantom{\alpha_{\iota(u)}}}) \coloneqq \ABFact{\(\tilde u\)}\)\;
  \uIf(\tcp*[f]{2nd condition for \(u=\tx^i, \tilde u=\overline\tx^j\)}){\(\iota(u) \ne \iota(\tilde u) \vee \alphabet(u) \ne \alphabet(\tilde u)\)}{
    \Return \(\min(\iota(u), \iota(\tilde u))\)\;
  }\uElseIf{\(\beta_1 = \tilde\beta_1 \wedge \cdots \wedge \beta_{\iota(u)} = \tilde\beta_{\iota(\tilde u)}\)}{
    \For(\tcp*[f]{solve \textsc{MaxSimK} for unary \(\alpha\) pairs}){\(i\in[\iota(u)]_0\)}{
      \(e_i \coloneqq\;\)\lIf{\(|\alpha_i| = |\tilde\alpha_i|\)}{\(\infty\)~\textbf{else}~\(\min(|\alpha_i|,|\tilde\alpha_i|)\)}
    }
    \Return \(\iota(u) + \min\{ e_i \mid i\in[\iota(u)]_0 \}\)\;
  }\Else{
    \Return \(\iota(u)\)\;
  }
\end{algorithm}

We can use \Cref{prop:binary:simon-characterization} to answer a number of questions regarding the structure of the congruence classes of \(\Sigma_2^* \mathord{/}{\sim_k}\).
For instance, for each \(w\) with \(|[w]_{\sim_k}| = \infty\), we have \(\tx^k \preceq w\) for some \(\tx\in\Sigma\) by the pigeonhole-principle.
The contrary is not true in general witnessed by the word \(v=\tb\tb\ta\tb\tb\) with respect to \(\sim_4\).
Its scattered factors of length four are \(\tb\tb\ta\tb, \tb\ta\tb\tb\) and \(\tb\tb\tb\tb\).
Therefore, each word in its class contains exactly one \(\ta\) (\(\ta\ta\not\preceq v\) but \(\ta\preceq v\) is), at least two \(\tb\) succeeding and preceding the \(\ta\) (\(\tb\tb\ta, \ta\tb\tb\preceq v\)) but not more than two \(\tb\) (\(\tb\tb\tb\ta,\ta\tb\tb\tb\not\preceq v\)).
Therefore, \(\tb\tb\ta\tb\tb\) is the only word in this class, but it contains \(\tb^4\).
By a famous result of Simon~\cite[Corollary~6.2.8]{lothaire1997:subwords}, all congruence classes of \(\sim_k\) are either infinite or singletons.
In the binary case, we can give a straightforward characterization of the finite/singleton and infinite classes.

\begin{theorem}\label{singleton}
  Let \(w \in \Sigma_2^*\), then \(|[w]_{\sim_k}| < \infty\). In particular, we have \(|[w]_{\sim_k}| = 1\) iff \(\iota(w) < k\) and \(|\alpha_i| < k - \iota(w)\) for all \(i\in[\iota(w)]_0\).
\end{theorem}

\ifpaper
\else
\begin{proof}
  Let \(w \in \Sigma_2^*\) with \(\iota(w) < k\), \(|\alpha_i| \leq k - \iota(w)\) for all \(i\in[\iota(w)]_0\), and \(\tilde w \in [w]_{\sim_k}\). By \Cref{rem:simon-refines-arch-num-equiv}, we have $\iota(w)=\iota(\tilde{w})=m$.
  By \Cref{prop:binary:simon-characterization}, we get \(\beta_i = \tilde \beta_i\) for all \(i\in[\iota(w)]\)
  and \([\alpha_i]_{\sim_{k-m}} = [\tilde \alpha_i]_{\sim_{k-m}}\) for all $i\in[\iota(w)]_0$.
  Additionally, we have by \Cref{prop:binary:beta-alpha-cases} and the fact that for \(u,v\in\Sigma_1^*\) we have \(u \sim_\ell v\), if and only if, \(\min(|u|, \ell) = \min(|v|, \ell)\) exactly \([\alpha_i]_{\sim_{k-m}} = [\tilde \alpha_i]_{\sim_{k-m}} = \{\tx^{|\alpha_i|}\}\) for all \(i\in[\iota(w)]_0\).
  Thus, $|[w]_{\sim_k}|=1$.

  Now let \(w \in \Sigma_2^*\) such that $\iota(w)\geq k$ or $|\alpha_i|\geq k-\iota(w)$ for some $i\in[\iota(w)]_0$.
  If \(\iota(w) \geq k\), then \(\Sigma_2^{\leq k} = \ScatFact_{\leq k}(w) \subseteq \ScatFact_{\leq k}(wv)\) and thus $wv\in[w]_{\sim_k}$ for all $v\in\Sigma_2^{\ast}$. This implies $|[w]_{\sim_k}|\neq 1$.
  On the other hand, if \(\iota(w) < k\) and there exists \(i\in[\iota(w)]_0\) with \(|\alpha_i| \geq k - \iota(w) \geq 1\), then \(\alpha_i \sim_{ k -\iota(w)} \alpha_i^j\) for all \(j \in \mathbb N\).
  Thus, we have \(\alpha_0\beta_1\alpha_1 \cdots \beta_i\alpha_i^j\beta_{i+1}\cdots \alpha_{\iota(w)-1}\beta_{\iota(w)}\alpha_{\iota(w)} \sim_k w\) for all \(j\in\N\) and we have again $|[w]_{\sim_k}|\neq 1$. \qed

\end{proof}

\fi


In the following, we will use \Cref{prop:binary:simon-characterization} to derive a formula for the precise value of \(|\Sigma_2^*\mathord{/}{\sim_k}|\).
Note that in the unary case, we have \(|\Sigma_1^*\mathord{/}{\sim_k}| = k + 1\) because the empty word has its own class.
By \Cref{rem:simon-refines-arch-num-equiv}, we know that there exists exactly one class w.r.t. $\sim_k$ of words with \(k\) arches.
We can consider the other classes by the common number of arches of their elements.
By \Cref{prop:binary:simon-characterization}, we can count classes based on the valid combinations of \(\beta\)-factors and number of classes for each \(\alpha\)-factors.
Because the \(\alpha\) are unary, we already know their number of classes.
These valid combinations are exactly given by \Cref{prop:binary:beta-alpha-cases}.
The first values for the number of classes separated by the number of arches are given in \Cref{tab:binary:simon-index-per-classes}.

\begin{table}[tb]
  \centering
  \caption{Index of \(\sim_k\) restricted to binary words with a fixed number of arches}
  \resizebox{0.8\columnwidth}{!}{
  \begin{tabular}{ccrrrrrrrrr}\toprule
    \multicolumn{1}{c}{} & \multicolumn{1}{c}{} & \multicolumn{9}{c}{Number of Arches} \\ \cmidrule(lr){3-11}
    &       & 0        & 1   & 2     & 3     & 4      & 5      & 6     & 7 & \(m\) \\ \midrule
    \multirow[c]{7}{*}{\rotatebox[origin=tr]{90}{Scat Fact Length}}
    & 1     & 3        &  1  &       &       &        &        &       &   \\
    & 2     & 5        & 10  & 1     &       &        &        &       &   \\
    & 3     & 7        & 26  & 34    & 1     &        &        &       &   \\
    & 4     & 9        & 50  & 136   & 116   & 1      &        &       &   \\
    & 5     & 11       & 82  & 358   & 712   & 396    & 1      &       &   \\
    & 6     & 13       & 122 & 748   & 2 564 & 3 728  & 1 352  & 1     &   \\
    & 7     & 15       & 170 & 1 354 & 6 824 & 18 364 & 19 520 & 4 616 & 1 \\ \midrule
    & \(k\) & \(2k+1\) &     &       &       &        &        &       &   \\ \bottomrule
  \end{tabular}
  }
  \label{tab:binary:simon-index-per-classes}
\end{table}

\begin{theorem}\label{prop:binary:simon-index-formula-per-archcount}
  The number of congruence classes of \(\Sigma_2^*\mathord{/}{\sim_k}\) of words with \(m < k\) arches (the entries of \Cref{tab:binary:simon-index-per-classes}) is given by 
  \[
    \norm*{
    \begin{pmatrix}
            k - m & k - m & k - m \\ 
            1     & 2     & 1     \\
            k - m & k - m & k - m
    \end{pmatrix}^{\!\!m}
    \mkern-8mu\cdot
    \begin{pmatrix}
            k - m \\ 1 \\ k - m
    \end{pmatrix}
    }_1 = c_k^m
  \]
  where \(c_k^{-1} \coloneqq 1\), \(c_k^0 \coloneqq 2k + 1\), and \(c_k^m \coloneqq 2\cdot (k - m + 1) \cdot c_{k-1}^{m-1} - 2\cdot (k - m) \cdot c_{k-2}^{m-2}\)
  where \(\norm{\cdot}_1\) denotes the \(1\)-norm.
\end{theorem}

\ifpaper
\begin{proof}
  First, we show that the matrix representation produces the correct values, then we show the characterization as recurrence.
  Note that \(k - m\) is fixed on the diagonals of \Cref{tab:binary:simon-index-per-classes}.
  Therefore, increasing both, increases just the exponent of the matrix.
  We show that the first column is correct and then proceed by induction along the diagonals.
  Denote the above matrix by \(D_{k,m}\).

  Let \(k \in\mathbb N_0\) and \(w \in \Sigma^*\) with \(m \coloneqq \iota(w) < k\).
  For \(i\in[m]_0\), all elements \(v \in [w]_{\sim_k}\) have \(k - m\) congruent \(\alpha_i\) by \Cref{prop:binary:simon-characterization}.
  By definition, their alphabets are proper subsets of \(\Sigma_2\).
  Therefore, they are either empty or non-empty unary words consisting of just \(\ta\) or \(\tb\).
  We separate the choice of \(\alpha_i\) into these three cases.
  Let \(M_\varepsilon^\ell \coloneqq \{ [w] \in \Sigma_2^*\mathord{/}{\sim_{(k - m) + \ell}} \mid \iota(w) = \ell, \alpha_0 \sim_{k-m} \varepsilon \}\) and \(M^\ell_{\tx} \coloneqq \{ [w] \in \Sigma_2^*\mathord{/}{\sim_{(k - m) + \ell}} \mid \iota(w) = \ell, \alpha_0 \sim_{k-m} \tx^r,r\in\N \}\) for \(\tx\in\Sigma_2\) be sets of \(\ell+(k-m)\) congruence classes of words with \(\ell\) arches, separated by the alphabet of \(\alpha_0\).
  Denote by \(e_{k, m} \coloneqq (|M_\ta^0|, |M_\varepsilon^0|, |M_\tb^0|)^\intercal = (k - m, 1, k - m)^\intercal\) the number of classes for zero arches.
  We show \(\norm{D_{k,m}^\ell \cdot e_{k,m}}_1 = (|M_\ta^\ell|, |M_\varepsilon^\ell|, |M_\tb^\ell|)^\intercal\).
  There are four choices for \(\beta_i\) which are given by \Cref{prop:binary:beta-alpha-cases}.
  Each choice of \(\beta_{i+1}\) depends on the preceding \(\alpha_i\) and limits the choices for the succeeding \(\alpha_{i+1}\).
  These are given by \Cref{prop:binary:beta-alpha-cases}, and correspond to the entries of the matrix because for \(\ell \geq 1\) we have
  \begin{align*}
    M_\varepsilon^\ell &= \{ [w]_{\sim_{(k - m) + \ell}} \in M_\varepsilon^\ell \mid \tx\in\Sigma_2, \beta_1(w) = \overline\tx\tx, \alpha_1(w) \sim_{k-m} \varepsilon \} \\
    &\phantom{=} \quad\sqcup \{ [w]_{\sim_{(k - m) + \ell}} \in M_\varepsilon^\ell \mid \tx\in\Sigma_2, \beta_1(w) = \overline\tx\tx, \alpha_1(w) \sim_{k-m} \tx^r,r\in\N \} \\
    &\cong \{\ta\tb,\tb\ta\} \times M_\varepsilon^{\ell-1} \sqcup \{\ta\tb\} \times M_{\tb}^{\ell-1} \sqcup \{\tb\ta\} \times M_{\ta}^{\ell-1} \\
    M_{\tx}^\ell &= \{ [w]_{\sim_{(k - m) + \ell}} \in M_\tx^\ell \mid \beta_1(w) = \overline\tx \} \\
    &\phantom{=} \quad\sqcup \{ [w]_{\sim_{(k - m) + \ell}} \in M_\tx^\ell \mid \beta_1(w) = \tx\overline\tx, \alpha_1(w) \sim_{k-m} \overline\tx^r,r\in\N \} \\
    &\phantom{=} \quad\sqcup \{ [w]_{\sim_{(k - m) + \ell}} \in M_\tx^\ell \mid \beta_1(w) = \tx\overline\tx, \alpha_1(w) \sim_{k-m} \varepsilon \} \\
    &\cong [k-m] \times \paren[\big]{\{\overline\tx\} \times M^{\ell-1}_\tx \sqcup \{\tx\overline\tx\} \times M^{\ell-1}_{\overline\tx} \sqcup \{\tx\overline \tx\} \times M^{\ell-1}_\varepsilon }.
  \end{align*}
  Therefore, each multiplication with the matrix increases the number \(\ell\) of arches  by one.
  Thus, for \(m = \ell\) we have the desired value as \(M_\varepsilon^m\) and \(M_\tx^m\) are sets of \(k\) congruence classes with \(m\) arches.
  Therefore, \(\norm{D_{k,m}^m \cdot e_{k,m}}_1\) corresponds to the number of classes with respect to \(\sim_k\) of words with \(m\) arches.

  The equivalence of the two formulas is left to show.
  The characteristic polynomial of \(D_{k,m}\) is given by
  $\chi_{D_{k,m}} = \det(D_{k,m} - \lambda I) = -\lambda^3 + 2 \lambda^2 + 2(k - m)\lambda^2 - 2 (k - m)\lambda$.
  By the Cayley-Hamilton theorem \cite{hamilton1853lectures}, \(D_{k,m}\) is a root of its characteristic polynomial and thus satisfies the recurrence
  \begin{align*}
    D_{k,m}^{\ell+2}
    &= 2\cdot D_{k,m}^{\ell+1} + 2\cdot(k - m)\cdot D_{k,m}^{\ell+1} - 2\cdot (k - m)\cdot D_{k,m}^{\ell} \\
    &= 2\cdot (k - m + 1)\cdot D_{k,m}^{\ell+1} - 2\cdot (k - m)\cdot D_{k,m}^\ell
  \end{align*}
  for \(\ell\in\mathbb N\).
  Note that \(e_{k, m} = e_{k + \ell, m + \ell}\) for all \(\ell \in \mathbb N_0\).
  Now we conclude by induction that
  \begin{align*}
    &\norm[\big]{D_{k+2, m+2}^{m+2}\cdot e_{k+2,m+2}}_1  
    = \norm[\big]{D_{k, m}^{m+2}\cdot e_{k,m}}_1 \\
    &= \norm[\big]{(2\cdot (k - m + 1)\cdot D_{k, m}^{m+1} - 2\cdot (k - m)\cdot D_{k, m}^m)\cdot e_{k,m}}_1 \\
    &= 2\cdot (k - m + 1)\cdot\norm[\big]{D_{k, m}^{m+1}\cdot e_{k,m}}_1 - 2\cdot (k - m)\cdot\norm[\big]{D_{k, m}^m\cdot e_{k,m}}_1 \\
    &= 2\cdot (k - m + 1)\cdot c_{k+1}^{m+1} - 2\cdot (k - m)\cdot c_{k}^{m} 
    = c_{k+2}^{m+2},
  \end{align*}
  because \(\norm{u \pm v}_1 = \norm{u}_1 \pm \norm{v}_1\) for all \(u = (u_i),v = (v_i)\in\mathbb R^n\) for which \(u_jv_j \geq 0\) for all \(j\in[n]\). \qed
\end{proof}

\else

\fi

\begin{remark}
  Note that by setting \(\Delta \coloneqq k - m\), the family of recurrences depends only on one variable \(\Delta\), because \(k - m = (k - \ell) - (m - \ell)\) holds for all $\ell\in\N$.
\end{remark}

\begin{remark}
	Some sequences in \Cref{tab:binary:simon-index-per-classes} are known sequences.
	The first and second diagonal are \href{https://oeis.org/A007052}{A007052} and \href{https://oeis.org/A018903}{A018903} resp. in \cite{oeis}.
	Both sequences are investigated in the work of Janjic~\cite{Janjic2010GeneralizedCW}.
	There, the two sequences appear as the number of compositions of \(n\in\mathbb N\), considering three (resp. five) differently colored \(1\)s.
	Furthermore, the sequences \(c_k^m\) seem to be equivalent to the family of sequences \((s_n)\) where \(s_0 = 1\) and \(s_1\) is fixed and \(s_{n+2}\) is the smallest number such that \(\frac{s_{n+2}}{s_{n+1}} > \frac{s_{n+1}}{s_n}\).
	These sequences where studied by Boyd~\cite{boyd1993pisot}.
\end{remark}

By \Cref{rem:simon-refines-arch-num-equiv}, we can count the number of classes separated by the universality of words with less than \(k\) arches.
This leads to the following immediate corollary which allows us to efficently calculate \(|\Sigma_2^*\mathord{/}{\sim_k}|\).

\begin{corollary}\label{prop:binary:simon-index}
  Let \(k\in \mathbb N_0\).
  Over a binary alphabet, the number of congruence classes of \(\sim_k\) is given by \(|\Sigma_2^*\mathord{/}{\sim_k}| = 1 + \sum_{m = 0}^{k-1} c_k^m\).
\end{corollary}
The first values of the sequence, some of which are already given in \cite{DBLP:journals/ipl/KarandikarKS15}, are
\begin{multline*}
                 1,
                 4,
                16,
                68,
               312,
             1 560,
             8 528,
            50 864,
           329 248,
         2 298 592,
        17 203 264,
       137 289 920, \\
     1 162 805 376,
    10 409 679 744,
    98 146 601 216,
   971 532 333 824,
10 068 845 515 264, \dots
\end{multline*}

\begin{table}[tb]
  \centering
  \caption[Number of classes of perfect universal binary words by number of arches]{Number of classes of perfect universal binary words restricted to a fixed number of arches}%
  \resizebox{0.8\columnwidth}{!}{
  \begin{tabular}{ccrrrrrrrrrr}\toprule
    \multicolumn{1}{c}{} & \multicolumn{1}{c}{} & \multicolumn{10}{c}{Number of Arches} \\ \cmidrule(lr){3-12}
    &       & 0        & 1   & 2     & 3     & 4      & 5      & 6     & 7 & 8 & \(m\) \\ \midrule
    \multirow[c]{7}{*}{\rotatebox[origin=tr]{90}{Scat Fact Length}}
    & 2     & 1     & 4            & 1   &       &       &        &        &       &   \\
    & 3     & 1     & 6            & 14  & 1     &       &        &        &       &   \\
    & 4     & 1     & 8            & 32  & 48    & 1     &        &        &       &   \\
    & 5     & 1     & 10           & 58  & 168   & 164   & 1      &        &       &   \\
    & 6     & 1     & 12           & 92  & 416   & 880   & 560    & 1      &       &   \\
    & 7     & 1     & 14           & 134 & 840   & 2 980 & 4 608  & 1 912  & 1     &   \\
    & 8     & 1     & 16           & 184 & 1 488 & 7 664 & 21 344 & 24 344 & 6 528 & 1 \\ \midrule
    & \(k\) & \(1\) & \(2k\)       &     &       &       &        &        &       &   \\ \bottomrule
  \end{tabular}
  }
  \label{tab:binary:perfect-universal-per-arch}
\end{table}

We can use the idea of \Cref{prop:binary:simon-index-formula-per-archcount} to count the number of perfect $k$-universal words, i.e., $k$-universal words with an empty rest (cf. \cite{fleischmann2021scattered}).
We can count them by replacing the vector from \Cref{prop:binary:simon-index-formula-per-archcount} with the initial distribution of \(\alpha_i\) values with \((0, 1, 0)^\intercal\).
Thus, the formula counts words starting or ending with an empty \(\alpha\).
Because the matrix does not change, we obtain the same recurrence with different initial values.
The \th{\(k\)} diagonal, shifted by one, is now given by the Lucas sequence of the first kind \(U(2\cdot k + 2, 2\cdot k)\), where \(U_n(P, Q)\) is given by
$U_0(P, Q) = 0$, $U_1(P, Q) = 1$, $U_n(P, Q) = P\cdot U_{n-1}(P, Q) - Q\cdot U_{n-2}(P, Q)$.
The first calculated values are given in \Cref{tab:binary:perfect-universal-per-arch}.
The first three diagonals of the table are the known integer sequences \href{https://oeis.org/A007070}{A007070}, \href{https://oeis.org/A084326}{A084326}, and \href{https://oeis.org/A190978}{A190978} in \cite{oeis}.

\section{Towards The Ternary Case of Simon's Congruence}\label{sec:ternary}
In the following, we will consider cases for the ternary alphabet based on the alphabets of the \(\alpha\)-factors with the goal of proving similar results to \Cref{prop:binary:beta-alpha-cases,prop:binary:simon-characterization} for ternary words, leading to \Cref{prop:ternary:beta-alpha-cases}.
By \Cref{prop:equiv-iff-all-aba-equiv}, it suffices to consider \(\alpha\beta\alpha\)-factors for characterizing congruence classes.
In \Cref{sec:alpha-beta} we already considered some cases for \(\alpha\beta\alpha\)-factors for arbitrary alphabets.
Note that if \(\m_1(w) = \w_1(w)\) then \(\core_1 = \varepsilon\).
Otherwise, if \(\m_1(w) \ne \w_1(w)\), then \(\core_1 \in (\Sigma \setminus \{\m_1(w), \w_1(w)\})^*\).
Thus, cores of ternary words are unary, and we denote the well-defined letter of the core by \(\ty\in\Sigma_3\).

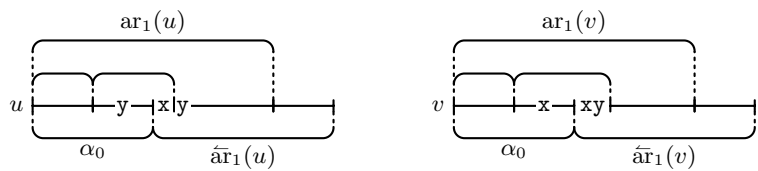
\begin{figure}
\centering
  \begin{tikzpicture}[scale=0.8]
    \draw[word] (0, 0) -- ++(1, 0) node [outerLabelStart] {\(u\)};
    \draw[word] (1, 0) -- ++(1, 0) node[innerLabel,yshift=-0.035cm] {\(\ty\)};
    \draw[word] (2.35, 0) -- ++(1.65, 0) node[pos=0,left,xshift=0.07cm] {\(\tx\)} node[innerLabel,pos=0,right,yshift=-0.035cm] {\(\ty\)};
    \draw[word] (4, 0) -- ++(1, 0);
    \draw[uword] (0, .9) -- ++(4, 0) node  {\(\ar_1(u)\)};
    \draw[dword] (2, -.35) -- ++(3, 0) node  {\(\ra_1(u)\)};
    \draw[dword] (0, -.35) -- ++(2, 0) node  {\(\alpha_0\)};

    \draw[uword] (0, .35) -- ++(1, 0);
    \draw[uword] (1, .35) -- ++(1.35, 0);

    \draw[word] (7, 0) -- ++(1, 0) node [outerLabelStart] {\(v\)};
    \draw[word] (8, 0) -- ++(1, 0) node[innerLabel,yshift=-0.00cm] {\(\tx\)};
    \draw[word] (9.6, 0) -- ++(1.4, 0) node[pos=0,left,yshift=-0.035cm,xshift=0.055cm] {\(\tx\ty\)};
    \draw[word] (11, 0) -- ++(1, 0);
    \draw[uword] (7, .9) -- ++(4, 0) node {\(\ar_1(v)\)};
    \draw[dword] (7, -.35) -- ++(2, 0) node  {\(\alpha_0\)};
    \draw[dword] (9, -.35) -- ++(3, 0) node {\(\ra_1(v)\)};
    \draw[uword] (7, .35) -- ++(1, 0);
    \draw[uword] (8, .35) -- ++(1.6, 0);

    \foreach \x in {1,2.35,8,9.6} {\draw[point] (\x, 0.1) -- ++(0, 0.2);}
    \foreach \x in {0,2,5,7,9,12} {\draw[point] (\x, -.1) -- ++(0, -.2);}
    \foreach \x in {4,11} {\draw[point] (\x, 0.1) -- ++(0, 0.8);}
    \foreach \x in {0,7} {\draw[point] (\x, 0.1) -- ++(0, 0.8);}
  \end{tikzpicture}
  \caption{Factorization of $\alpha$ in the ternary case assuming \(\core\in\ty^+\).}
  \label{fig:alpha-small-arches}
\end{figure}

We use a variant of the  {\em Kronecker-$\delta$} for a boolean predicate $P$ as $\delta_{P(x)}=1$ if $P(x)$ is true and $0$ otherwise to express a condition on the alphabet of the rest of a binary \(\alpha\)-factor (cf.~\Cref{fig:alpha-small-arches}).
If an \(\alpha_0\)'s rest contains the letter \(\ty\) different from the reverse modus \(\tx \coloneqq \w(w)\), then $\re(\alpha_0)\w(w)$ builds another arch ending before the core (left).
This lowers the level of congruence, up to which we can determine the \(\core\), by one.
If \(\ty \npreceq \re(\alpha_0)\) the next \(\ty\) is in the \(\core\) (right).

We always assume that \(k \geq 2\) because we characterize the congruence of \(1\)-universal words. Moreover, let $w,\tilde{w}\in\Sigma_3^{\ast}$ with $1=\iota(w)=\iota(\tilde{w})$.

First, we prove a useful lemma which characterizes the congruence of two ternary words with the same modus and reverse modus. 
Together with \Cref{prop:different-modi-force-alpha-univ-k-1}, this immediately implies several cases.
\begin{lemma}\label{prop:ternary:same-modi-imply-equiv-iff-core-equiv}
  Let \(\m(w) = \m(\tilde w)\) and \(\w(w) = \w(\tilde w)\), we have \(w \sim_k \tilde w\) iff \(\alpha_i \sim_{k-1} \tilde \alpha_i\) for all \(i\in[1]_0\) and \(\core_1 \sim_{k-c} \widetilde\core_1 \in \ty^*\) where \(c \coloneqq \iota(\alpha_0) + \delta_{\ty\preceq\re(\alpha_0)} + \iota(\alpha_1) + \delta_{\ty\preceq\er(\alpha_1)}\).
\end{lemma}
\ifpaper
\else
\begin{proof}
  Assume without loss of generality \(\beta = 01^\ell2\) and \(\tilde\beta = 01^{\tilde \ell}2\) for some \(\ell,\tilde \ell\in\mathbb N_0\) and \(\ell \leq \tilde \ell\).

  Further, assume \(\alpha_i \sim_{k-1} \tilde\alpha_i\) for all \(i\in[1]_0\) and \(\min(k - c, \ell) = \min(k - c,\tilde\ell)\) as above.
  By \Cref{prop:simon-equiv-iff-equiv-with-same-alpha}, we can assume that \(\alpha_i = \tilde \alpha_i\) for all \(i\in[1]_0\).
  If \(\ell = \tilde \ell\) or if \(\ell < k - c\) and we thus have \(\ell = \tilde \ell\), we have \(w = \tilde w\).
  Hence, assume \(k - c \leq \ell < \tilde \ell\).
  We show that \(\alpha_0 \cdot 01^\ell2 \cdot \alpha_1 \sim_k \alpha_0 \cdot 01^{\ell+1}2 \cdot \alpha_1\) using \Cref{prop:simon:letter-elim}, then the claim follows by repeated application of the lemma.
  Note that \(k - \ell \leq c = \iota_{\{0,1\}}(\alpha_00) + \iota_{\{1,2\}}(2\alpha_1)\).
  Therefore, we can factorize \(\alpha_00\) and \(2\alpha_1\) into at least \(c\) factors with binary alphabet each.
 The claim follows with \(c + \ell \geq k\).

  Because \(\alpha_i \sim_{k-1} \tilde\alpha_i\) is a necessary condition for \(w \sim_k \tilde w\) by \Cref{prop:simon-equiv-implies-equiv-aba-factors}, assume this condition holds and \(\min(k - c, \ell) \ne \min(k - c,\tilde\ell)\).
  Again, by \Cref{prop:simon-equiv-iff-equiv-with-same-alpha} assume that \(\alpha_i = \tilde \alpha_i\) for all \(i\in[1]_0\).
  Thus, \(\ell < k - c\) and \(\ell < \tilde \ell\).
  Define \(v \coloneqq \m_{\{0,1\}}(\alpha_00)\cdot 1^{\min(k-c,\tilde\ell)} \cdot \w_{\{1,2\}}(2\alpha_1) \in \Sigma_3^{\leq k}\).
  By construction, we have \(v \npreceq w\) because \(\ell < \min(k-c,\tilde\ell)\) and \(v \preceq \tilde w\) because \(\tilde \ell \geq \min(k-c,\tilde\ell)\).
  Thus, \(w \nsim_k \tilde w\). \qed
\end{proof}

\fi

Since in the ternary case, there are congruent words having different modi or reverse modi,
Lemma~\ref{prop:ternary:same-modi-imply-equiv-iff-core-equiv} does not imply a full characterization.
This leads to two cases in the following classification (case 3 and 5 out of the 9 cases in \Cref{tab:ternary:alpha-beta-cases}).
These two cases correspond to the first case in the following theorem.

\begin{theorem}\label{prop:ternary:beta-alpha-cases}
  For \(w, \tilde w \in \Sigma_3^*\) we have \(w\sim_k \tilde w\) iff \(\alpha_i\sim_{k-1} \tilde\alpha_i\) for all \(i\in[1]_0\), and one of the following\\
  1. \(|\alphabet(\alpha_i)| = 2\), \(\alphabet(\alpha_{1-i}) \cap \alphabet(\alpha_i) = \emptyset\), and \(\iota(\alpha_i) \geq k-1\) for some \(i\in[1]_0\),\\
  2. \(\m(w) = \m(\tilde w)\), \(\w(w) = \w(\tilde w)\), \(\core \sim_{k - c} \widetilde\core\) with \(c \coloneqq \iota(\alpha_0) + \delta_{\ty\preceq \alpha_0} + \iota(\alpha_1) + \delta_{\ty \in \alpha_1}\).\\
  For all possibilities distinguishing the \(\beta\)-factors, see \Cref{tab:ternary:alpha-beta-cases}.
\end{theorem}
\ifpaper
\else
\begin{proof}
  Above, we already considered the cases \(\alpha_0 = \alpha_1 = \varepsilon\) (\Cref{prop:both-alpha-empty}) and \(\letters(\alpha_0) = \letters(\alpha_1) \in \binom{\Sigma}{|\Sigma|-1}\) (\Cref{prop:alpha-miss-same-single}).
  Next we cover the case of the two alphabets of size two intersecting, finishing the \(2 = |\letters(\alpha_0)| = |\letters(\alpha_1)|\) case.
  Notice that in the case where both \(\alpha\)-factors are binary, we only have \(|\letters(\alpha_0) \cap \letters(\alpha_1)| = 1\) left to consider.
  
  \begin{numclaim}\label{prop:ternary:both-alpha-miss-diff-single}
    Let \(|\letters(\alpha_0)| = |\letters(\alpha_1)| = 2\).
    We have \(w \sim_k \tilde w\) iff \(\alpha_i \sim_{k-1} \tilde\alpha_i\) for all \(i\in[1]_0\) and \(\core_1 \sim_{k-c} \widetilde\core_1\), where \(c \coloneqq  \iota(\alpha_0) + \delta_{\ty\preceq \re(\alpha_0)} + \iota(\alpha_1) + \delta_{\ty \preceq\er(\alpha_1)}\).
  \end{numclaim}
  \begin{proof}
  Because \(|\alphabet(\alpha_i)| = 2\) and in both directions we have \(\letters(\alpha_i) = \letters(\tilde \alpha_i)\) for all \(i\in[1]_0\), the claim follows by \Cref{prop:ternary:same-modi-imply-equiv-iff-core-equiv}. \qed
\end{proof}

  The case described in \Cref{prop:ternary:both-alpha-miss-diff-single} is problematic for the enumeration of classes, because the choices of \(\beta\) depend on the universality and the rests' alphabets of the surrounding \(\alpha\)-factors.
  This is in contrast to the binary case, where we only had to distinguish between empty and the two types of unary \(\alpha\)-factors.
  
  Next, we consider the general case \(2 = |\letters(\alpha_0)|, 1 = |\letters(\alpha_1)|\).
  Due to symmetry, it suffices to consider only the cases of the form \(|\letters(\alpha_0)| \geq |\letters(\alpha_1)|\).
  The first cases is similar to the one of  \Cref{prop:ternary:both-alpha-miss-diff-single} because the modi are fixed by the structure of the \(\alpha\)-factors.
  
  \begin{numclaim}\label{prop:ternary:alpha-miss-two-one-overlapping}
   Let \(|\alphabet(\alpha_0)| = 2\), \(|\alphabet(\alpha_1)| = 1\) and \(\alphabet(\alpha_0) \cap \alphabet(\alpha_1) \ne \emptyset\).
    Then, \(w \sim_k \tilde w\) iff  \(k=1\) or \(\alpha_i \sim_{k-1} \tilde\alpha_i\) for all \(i\in[1]_0\) and \(\core_1 \sim_{k-c} \widetilde\core_1\), where \(c \coloneqq  \iota(\alpha_0) + \delta_{\ty\preceq \re(\alpha_0)} + \delta_{\ty \preceq\er(\alpha_1)}\).
  \end{numclaim}
  
  It is left to consider the case where the two binary alphabets intersect.
  
  \begin{numclaim}\label{prop:ternary:alpha-miss-two-one-disjoint}
    Let  \(|\alphabet(\alpha_0)| = 2\), \(|\alphabet(\alpha_1)| = 1\) and \(\alphabet(\alpha_0) \cap \alphabet(\alpha_1) = \emptyset\).
    Then, \(w \sim_k \tilde w\) iff either
     \(\iota(\alpha_0) = k - 1 = \iota(\tilde{\alpha}_0)\), \(\alpha_1 \sim_{k-1} \tilde\alpha_1\),
    or
    \(\alpha_i\sim_{k-1}\tilde\alpha_i\) for all \(i\in[1]_0\), \(\iota(\alpha_0) < k - 1\), \(\w(w) = \w(\tilde w)\), \(\core_1 \sim_{k-c} \widetilde\core_1 \in \ty^+\) where \(c \coloneqq \iota(\alpha_0) + \delta_{\ty \preceq \er(\alpha_0)}\).
  \end{numclaim}
  \begin{proof}
  Let \(w, \tilde w\in\Sigma_3^*\) as above and assume \(w \sim_k \tilde w\) and without loss of generality that \(\alpha_0 \in \{0,1\}^*\), \(\m(w) = \m(\tilde w) = 2\) and \(\alpha_1\in 2^+\).
  By \Cref{prop:simon-equiv-implies-equiv-aba-factors}, we have \(\alpha_i\sim_{k-1}\tilde\alpha_i\) for all \(i\in[1]_0\).
  If \(\iota(\alpha_1) \geq k - 1\), then we are done.
  Thus, assume \(\iota(\alpha_1) < k - 1\).
  We cannot have \(\w(w) \ne \w(\tilde w)\) by \Cref{prop:different-modi-force-alpha-univ-k-1}.
  Thus, we have \(\w(w) = \w(\tilde w)\) as well as \(\m(w) = \m(\tilde w)\) and the claim follows by \Cref{prop:ternary:same-modi-imply-equiv-iff-core-equiv}.

  Now, for the other direction assume \(\alpha_i\sim_{k-1}\tilde\alpha_i\) for all \(i\in[1]_0\).
  First, we consider \(\iota(\alpha_0) = k - 1 = \iota(\tilde\alpha_0)\).
  We have \(\beta = v2\) and \(\tilde \beta = \tilde v 2\) for \(v, \tilde v \in \{0, 1\}^*\).
  Because \(\alpha_0 \sim_{k-1} \tilde\alpha_0 \sim_{k-1} (01)^{k-1}\) we have \(\alpha_0 v \sim_k \alpha_0\tilde v \sim_k \tilde\alpha_0 \tilde v\) by \Cref{prop:extend-equiv-words-by-univ-word-increases-equiv} for \(\{0,1\}\).
  Therefore, we have \(w = \alpha_0 v 2\alpha_1 \sim_k \tilde \alpha_0 \tilde v 2 \tilde\alpha_1 = \tilde w\) because \(\sim_k\) is a congruence.
  The other case follows again from \Cref{prop:ternary:same-modi-imply-equiv-iff-core-equiv}. \qed
\end{proof}

  The next lemma examines the case \(2 = |\alphabet(\alpha_0)|, 0 = |\alphabet(\alpha_1)|\).
  This case is similar to the one above.
  If \(\alpha_1\) has at least \(k-1\) arches, then we can permute the modi.
  Otherwise, everything is fixed and the level of congruence of the cores is determined by the structure of the \(\alpha\)-factors. The proof is analogous to the one of \Cref{prop:ternary:alpha-miss-two-one-disjoint} because in the proof we do not use the fact that \(\alpha_0\) is non-empty.
  
  \begin{numclaim}\label{prop:ternary:alpha-binary-empty}
    Let  \(|\alphabet(\alpha_0)| = 2\), and \(\alpha_1 = \varepsilon\).
    Then, \(w \sim_k \tilde w\) iff either
     \(\iota(\alpha_0) = k - 1 = \iota(\tilde{\alpha}_0)\),
    or
    \(\alpha_0\sim_{k-1}\tilde\alpha_0\), \(\iota(\alpha_0) < k - 1\), \(\w(w) = \w(\tilde w)\), and \(\core_1 \sim_{k-c} \widetilde\core_1 \in \ty^+\) where \(c \coloneqq \iota(\alpha_0) + \delta_{\ty \preceq\er(\alpha_0)}\).
  \end{numclaim}
  
  With \Cref{prop:ternary:alpha-binary-empty} the case distinction for one of the $\alpha$ containing two letters, is completed. In the following cases we investigate the situation, where both \(\alpha\) have maximal 1-letter alphabets.
  
  \begin{numclaim}\label{prop:ternary:alpha-unary-empty}
    Let  \(|\alphabet(\alpha_0)| = 1\), and \(\alpha_1 = \varepsilon\).
    Then, \(w \sim_k \tilde w\) iff \(\alpha_1 \sim_{k-1} \tilde \alpha_1\), \(\m(w) = \m(\tilde w)\), \(\w(w) = \w(\tilde w)\) and \(\core_1 \sim_{k-c} \widetilde\core_1 \in \ty^+\) where \(c \coloneqq \delta_{\ty\preceq \alpha_1}\).
  \end{numclaim}
  \begin{proof}
  Assume \(w \sim_k \tilde w\).
  By \Cref{prop:different-modi-force-alpha-univ-k-1}, we can assume that the modi are the same.
  Now the claim, as well as, the other direction follow from \Cref{prop:ternary:same-modi-imply-equiv-iff-core-equiv}. \qed
\end{proof}

  The last two cases follow analogously to \Cref{prop:ternary:alpha-unary-empty}.
  
  \begin{numclaim}\label{prop:ternary:alpha-both-unary-overlapping}
    Let  \(\alphabet(\alpha_1) = \alphabet(\alpha_0) \in \Sigma_3\).
    Then, \(w \sim_k \tilde w\) iff \(\alpha_1 \sim_{k-1} \tilde \alpha_1\), \(\m(w) = \m(\tilde w)\), \(\w(w) = \w(\tilde w)\) and \(\core_1 \sim_{k-c} \widetilde\core_1 \in \ty^+\) where \(c \coloneqq \delta_{\ty\preceq \alpha_0} + \delta_{\ty \in \alpha_1}\).
  \end{numclaim}
  
  \begin{numclaim}\label{prop:ternary:alpha-both-unary-disjoint}
    Let \(\alphabet(\alpha_1) \ne \alphabet(\alpha_0)\) and $|\letters(\alpha_1)|=|\letters(\alpha_0)|=1$.
    Then, \(w \sim_k \tilde w\) iff \(\alpha_1 \sim_{k-1} \tilde \alpha_1\), \(\m(w) = \m(\tilde w)\), \(\w(w) = \w(\tilde w)\) and \(\core_1 \sim_{k-c} \widetilde\core_1 \in \ty^+\) where \(c \coloneqq \delta_{\ty\preceq \alpha_0} + \delta_{\ty \in \alpha_1}\).
  \end{numclaim}
  This concludes the proof. \qed
\end{proof}
\fi

\begin{table}
  \caption{The possibilities for the $\beta$-factor of $w=\alpha_0\beta\alpha_1$, assuming $\ta,\tb,\tc\in\Sigma_3$ different.
  Note that in the cases (1,1), (1,0), (0,0) the letters not fixed by the $\alpha$-factors can be chosen arbitrarily but differently from $\Sigma_3$.}
  \begin{center}
  \resizebox{0.8\columnwidth}{!}{
      \begin{tabular}{cllll} \toprule
        \emph{\(\abs{\alphabet(\alpha_0)}\),\(\abs{\alphabet(\alpha_1)}\)}                      & \(\letters(\alpha_0)\) & \(\letters(\alpha_1)\) & \(\beta\) {RegExp}                       & {Stated In} \\ \midrule
        \multirow[c]{2}{*}{\(2\),\(2\)} & $\{\ta,\tb\}$          & $\{\ta,\tc\}$          & $\tb\ta^{\ast}\tc$                              & {Prop.~\ref{prop:alpha-miss-same-single}} \\
                                         & $\{\ta,\tb\}$          & $\{\ta,\tb\}$          & $\tc$                                           &  \\ \cmidrule(lr){1-5}
        \multirow[c]{2}{*}{\(2\),\(1\)} & $\{\ta,\tb\}$          & $\{\tc\}$& $(\ta\tb^+\mid \tb\ta^+)\tc$                    &  \\
                                         & $\{\ta,\tb\}$          & $\{\ta\}$          & $\tb\ta^{\ast}\tc$                              &  \\ \cmidrule(lr){1-5}
        \(2\),\(0\)                     & $\{\ta,\tb\}$          & $\emptyset$          & $(\ta\tb^+ \mid \tb^+\ta)\tc$                    & \\ \cmidrule(lr){1-5} 
        \multirow[c]{2}{*}{\(1\),\(1\)} & $\{\ta\}$              & $\{\tb\}$              & $\ta\tb^+\tc \mid \ta\tc^+\tb \mid \tc\ta^+\tb$ & \\
                                         & $\{\ta\}$              & $\{\ta\}$              & $\tb \ta^*\tc$                                  & \\ \cmidrule(lr){1-5} 
        \(1\),\(0\)                     & $\{\ta\}$              & $\emptyset$           & $\tb\ta^*\tc \mid \ta\tb^+\tc $                 & \\ \cmidrule(lr){1-5} 
        \(0\),\(0\)                     & $\emptyset$            & $\emptyset$            & $\ta\tb^+\tc $                                  & {Prop.~\ref{prop:both-alpha-empty}} \\ \bottomrule
      \end{tabular}
      }
    \end{center}

  \label{tab:ternary:alpha-beta-cases}
\end{table}


\section{Conclusion}\label{sec:conclusion}

In 2021, Kosche et al.~\cite{kosche2021absent} first implicitly used a new factorization to find absent scattered factors in words algorithmically.
Later, in 2022, Fleischmann et al.~\cite{DBLP:conf/dcfs/FleischmannHHMN22} introduced this factorization as \(\alpha\)-\(\beta\)-factorization and used it to investigate the classes of Simon's congruence separated by the number of shortest absent scattered factors, to characterize the classes for arbitrary alphabets for some fixed numbers of shortest absent scattered factors and to give explicit formulas for these subsets.
In this paper, we investigated the \(\alpha\)-\(\beta\)-factorization as an object of intrinsic interest. This leads to a result characterizing \(k\)-congruence of \(m\)-universal words in terms of their \(1\)-universal \(\alpha\beta\alpha\)-factors.
In the case of the binary and ternary alphabet, we fully characterized the congruence of words in terms of their single factors.
Moreover, using this characterization, we gave a formula for the number of classes of binary words for each \(k\), characterized the finite classes, and gave a conceptually simple linear time algorithm for testing \textsc{MaxSimK} for binary words.

The modus of the layered arch factorizations used in the proof of \Cref{prop:different-modi-force-alpha-univ-k-1} and throughout the literature~\cite{DBLP:conf/automata/Simon75,lothaire1997:subwords,DBLP:conf/mfcs/FleischerK18}, can be regarded as the optimal word to jump to certain letters in certain parts of the word.
The \(\alpha\)-\(\beta\)-factorization encapsulates the first layer (arches w.r.t. \(\Sigma\)) of these factorizations for all indicies.
For small alphabets (this paper) and shortest abscent scattered factors (c.f.~\cite{DBLP:conf/dcfs/FleischmannHHMN22}) this allows the characterization and enumeration of classes.
Extending this idea to lower layers (arches w.r.t. some \(\Omega \subset \Sigma\)), is left as future work.

\bibliographystyle{splncs04}
\bibliography{references}

\begin{thebibliography}{10}
\providecommand{\url}[1]{\texttt{#1}}
\providecommand{\urlprefix}{URL }
\providecommand{\doi}[1]{https://doi.org/#1}

\bibitem{DBLP:conf/dlt/BarkerFHMN20}
Barker, L., Fleischmann, P., Harwardt, K., Manea, F., Nowotka, D.: {S}cattered
  {F}actor-{U}niversality of {W}ords. In: DLT (2020)

\bibitem{boyd1993pisot}
Boyd, D.W.: Linear {R}ecurrence {R}elations for some {G}eneralized {P}isot
  {S}equences. In: {C}{N}{T}{A}. pp. 333--340 (1993)

\bibitem{day2021edit}
Day, J., Fleischmann, P., Kosche, M., Ko{\ss}, T., Manea, F., Siemer, S.: The
  {E}dit {D}istance to k-{S}ubsequence {U}niversality. In: {STACS}. vol.~187,
  pp. 25:1--25:19 (2021)

\bibitem{DBLP:conf/mfcs/FleischerK18}
Fleischer, L., Kufleitner, M.: Testing {S}imon's {C}ongruence. In: Potapov, I.,
  Spirakis, P.G., Worrell, J. (eds.) {MFCS}. vol.~117, pp. 62:1--62:13 (2018)

\bibitem{fleischmann2021scattered}
Fleischmann, P., Germann, S., Nowotka, D.: Scattered {F}actor {U}niversality --
  {T}he {P}ower of the {R}emainder. preprint arXiv:2104.09063 (published at
  RuFiDim)  (2021)

\bibitem{DBLP:conf/dcfs/FleischmannHHMN22}
Fleischmann, P., Haschke, L., Huch, A., Mayrock, A., Nowotka, D.: Nearly
  {$k$}-{U}niversal {W}ords -- {I}nvestigating a {P}art of {S}imon's
  {C}ongruence. In: {DCFS}. {LNCS}, vol. 13439, pp. 57--71. Springer (2022)

\bibitem{hamilton1853lectures}
Frobenius, H.: {\"U}ber lineare Substitutionen und bilineare Formen, vol.~1878.
  De Gruyter (1878)

\bibitem{DBLP:conf/stacs/GawrychowskiKKM21}
Gawrychowski, P., Kosche, M., Ko{\ss}, T., Manea, F., Siemer, S.: Efficiently
  {T}esting {S}imon's {C}ongruence. In: {STACS}. LIPIcs, vol.~187, pp.
  34:1--34:18 (2021)

\bibitem{DBLP:journals/tcs/Hebrard91}
H{\'{e}}brard, J.: An {A}lgorithm for {D}istinguishing {E}fficiently
  {B}it-strings by their {S}ubsequences. TCS  \textbf{82}(1),  35--49 (1991)

\bibitem{Janjic2010GeneralizedCW}
Janjic, M.: Generalized {C}ompositions with a {F}ixed {N}umber of {P}arts.
  arXiv:1012.3892  (2010)

\bibitem{karandikar2016height}
Karandikar, P., Schnoebelen, P.: The {H}eight of {P}iecewise-{T}estable
  {L}anguages with {A}pplications in {L}ogical {C}omplexity. In: CSL (2016)

\bibitem{DBLP:journals/ipl/KarandikarKS15}
Karandikar, P., Kufleitner, M., Schnoebelen, P.: On the {I}ndex of {S}imon's
  {C}ongruence for {P}iecewise {T}estability. Information Processing Letters
  \textbf{115}(4),  515--519 (2015)

\bibitem{journals/lmcs/KarandikarS19}
Karandikar, P., Schnoebelen, P.: The {H}eight of {P}iecewise-{T}estable
  {L}anguages and the {C}omplexity of the {L}ogic of {S}ubwords. LMCS
  \textbf{15}(2) (2019)

\bibitem{DBLP:conf/dcfs/KimHKS22}
Kim, S., Han, Y., Ko, S., Salomaa, K.: On simon's congruence closure of a
  string. In: Han, Y., Vaszil, G. (eds.) {DCFS}. {LNCS}, vol. 13439, pp.
  127--141. Springer (2022)

\bibitem{DBLP:conf/isaac/KimKH22}
Kim, S., Ko, S., Han, Y.: Simon's congruence pattern matching. In: Bae, S.W.,
  Park, H. (eds.) ISAAC. LIPIcs, vol.~248, pp. 60:1--60:17. Schloss Dagstuhl -
  Leibniz-Zentrum f{\"{u}}r Informatik (2022)

\bibitem{kosche2021absent}
Kosche, M., Ko{\ss}, T., Manea, F., Siemer, S.: Absent {S}ubsequences in
  {W}ords. In: RP. pp. 115--131. Springer (2021)

\bibitem{pach2}
Kátai-Urbán, K., Pach, P., Pluhár, G., Pongrácz, A., Szabó, C.: On the
  word problem for syntactic monoids of piecewise testable languages. Semigroup
  Forum  \textbf{84}(2) (2012)

\bibitem{lothaire1997}
Lothaire, M.: Combinatorics on Words. Cambridge Mathematical Library, Cambridge
  University Press, 2 edn. (1997)

\bibitem{oeis}
{OEIS Foundation Inc.}: The {O}n-{L}ine {E}ncyclopedia of {I}nteger {S}equences
  (2022), \url{http://oeis.org}

\bibitem{pach}
Pach, P.: Normal forms under {S}imon's congruence. Semigroup Forum
  \textbf{97}(2) (2018)

\bibitem{lothaire1997:subwords}
Sakarovitch, J., Simon, I.: Subwords, chap.~6, pp. 105--144. In: Cambridge
  Mathematical Library \cite{lothaire1997}, 2 edn. (1997)

\bibitem{simon1972thesis}
Simon, I.: Hierarchies of Events with Dot-depth One. Ph.D. thesis, University
  of Waterloo, Department of Applied Analysis and Computer Science (1972)

\bibitem{DBLP:conf/automata/Simon75}
Simon, I.: Piecewise {T}estable {E}vents. In: Barkhage, H. (ed.) Automata
  Theory and Formal Languages, 2nd {GI} Conference. LNCS, vol.~33, pp.
  214--222. Springer (1975)

\end{thebibliography}

%

\end{document}